%% file: rdmd.tex
\documentclass[onefignum,onetabnum]{siamonline171218}

\input{ex_shared}

\ifpdf
\hypersetup{
  pdftitle={Randomized Dynamic Mode Decomposition},
  pdfauthor={Erichson, Mathelin, Kutz, Brunton}
}
\fi

\begin{document}

\maketitle

\begin{abstract}
This paper presents a randomized algorithm for computing the near-optimal low-rank dynamic mode decomposition (DMD). 
Randomized algorithms are emerging techniques to compute low-rank matrix approximations at a fraction of the cost of deterministic algorithms, easing the computational challenges arising in the area of `big data'. 
The idea is to derive a small matrix from the high-dimensional data, which is then used to efficiently compute the dynamic modes and eigenvalues. 
The algorithm is presented in a modular probabilistic framework, and the approximation quality can be controlled via oversampling and power iterations. 
The effectiveness of the resulting randomized DMD algorithm is demonstrated on several benchmark examples of increasing complexity, providing an accurate and efficient approach to extract spatiotemporal coherent structures from big data in a framework that scales with the intrinsic rank of the data, rather than the ambient measurement dimension. For this work we assume that the dynamics of the problem under consideration is evolving on a low-dimensional subspace that is well characterized by a fast decaying singular value spectrum.
\end{abstract}



\vspace{+1cm}
\section[Introduction]{Introduction}
Extracting dominant coherent structures and modal expansions from high-dimensional data is a cornerstone of computational science and engineering~\cite{Taira2017aiaa}.  
The dynamic mode decomposition (DMD) is a leading data-driven algorithm to extract spatiotemporal coherent structures from high-dimensional data sets~\cite{Schmid2010jfm,Tu2014jcd,Kutz2016book}. 
DMD originated in the fluid dynamics community~\cite{Schmid2010jfm,Rowley2009jfm}, where the identification of coherent structures, or \emph{modes}, is often an important step in building models for prediction, estimation, and control~\cite{Brunton2015amr,Taira2017aiaa}.  
In contrast to the classic proper orthogonal decomposition (POD)~\cite{Berkooz:1993,HLBR_turb}, which orders modes based on how much energy or variance of the flow they capture, DMD identifies spatially correlated modes that oscillate at a fixed frequency in time, possibly with an exponentially growing or decaying envelope.  
Thus, DMD combines the advantageous features of POD in space and the Fourier transform in time~\cite{Kutz2016book}.

With rapidly increasing volumes of measurement data from simulations and experiments, modal extraction algorithms such as DMD may become prohibitively expensive, especially for online or real-time analysis.  
Even though DMD is based on an efficient singular value decomposition (SVD), computations scale with the dimension of the measurements, rather than with the intrinsic dimension of the data.  
With increasingly vast measurements, it is often the case that the intrinsic rank of the data does not increase appreciably, even as the dimension of the ambient measurements grows.  
Indeed, many high-dimensional systems exhibit low-dimensional patterns and coherent structures, which is one of the driving perspectives in both POD and DMD analysis.

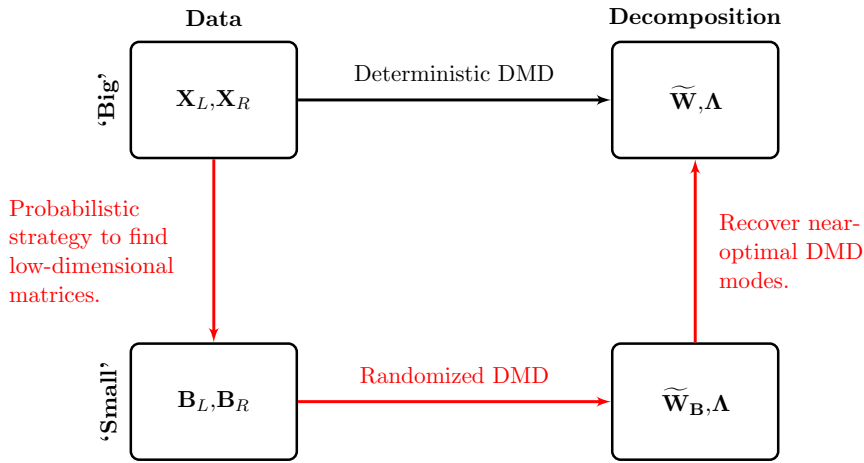
\begin{figure}[!b]
	\begin{center}
		\scalebox{0.78}{
			\begin{tikzpicture}[auto,node distance = 2cm,>=latex']
			\node [blocky4,name=fdata] { ${\bX}_L$,${\bX}_R$ };       
			\node [blocky4,name=fdmd, right of=fdata,node distance=8cm]  {${\bW}$,$\mathbf{\Lambda}$};
			\node [blocky4,name=cdata,below of=fdata,node distance=5.0cm] {${\bB}_L$,${\bB}_R$};
			\node [blocky4,name=cdmd,below of=fdmd,node distance=5.0cm] {$\widetilde{\bW}_\bB$,$\mathbf{\Lambda}$};
			\node [blocky5,name=dummyleft,below of=fdata,node distance=1.5cm] {};
			\node [blocky5,name=dummyright,below of=fdmd,node distance=1.5cm] {};
			
			\path [draw, ->, line width=1.4pt] (fdata) -- node[name=toparrow, above of=fdata, node distance=.45cm]{{Deterministic DMD}} (fdmd);
			\path [draw, ->, line width=1.4pt, red] (cdata) -- node[name=bottomarrow,above of=cdata,node distance=.45cm]{{Randomized DMD}}(cdmd);
			\path [draw, ->, line width=1.4pt, red] (fdata) -- node[name=leftarrow, left of=fdata, node distance = 1.9cm, text width=3.0cm]{Probabilistic strategy to find low-dimensional matrices.}(cdata);
			\path [draw, <-, line width=1.4pt, red] (fdmd) -- node[name=rightarrow, right of=fdmd, node distance = 1.9cm, text width=3.0cm]{Recover near-optimal DMD modes. }(cdmd);
			\node [name=data,above of=fdata,node distance=1.35cm] {\textbf{Data}};
			\node [name=modes,above of=fdmd,node distance=1.35cm] {\textbf{Decomposition}};
			\node [name=full,left of=fdata,node distance=1.75cm]{\begin{sideways}\textbf{`Big'}\end{sideways}};
			\node [name=full,left of=cdata,node distance=1.75cm]{\begin{sideways}\textbf{`Small'}\end{sideways}};
			
			\end{tikzpicture}}
	\end{center}
	\caption{Conceptual architecture of the randomized dynamic mode decomposition (rDMD). First, small matrices are derived from the high-dimensional input data. The low-dimensional snapshot matrices $\mathbf{B}_L$ and $\mathbf{B}_R$ are then used to compute the approximate dynamic modes $\widetilde{\bW}_\bB$, and eigenvalues $\mathbf{\Lambda}$. Finally, the near-optimal modes ${\bW}$ may be recovered.}
	\label{Fig:randLA}
\end{figure}

\newpage
In this work, we develop a computationally effective strategy to compute the DMD using randomized linear algebra, which scales with the intrinsic rank of the dynamics, rather than with the measurement dimension.
More concretely, we embed the DMD in a probabilistic framework by following the seminal work of Halko et al.~\cite{halko2011rand}. The main concept is depicted in Figure~\ref{Fig:randLA}. Numerical experiments show that our randomized algorithm achieves considerable speedups over previously proposed algorithms for computing the DMD such as the compressed DMD algorithm~\cite{Brunton2015jcd}. Further, the approximation error can be controlled via oversampling and additional power iterations. This allows the user to choose the optimal trade-off between computational time and accuracy. 
Importantly, we also demonstrate the ability to handle data which are too big to fit into fast memory by using a blocked matrix scheme. 

\subsection{Related Work}
\begin{figure}[!b]
	\centering
	\DeclareGraphicsExtensions{.pdf}
	\includegraphics[width=.6\textwidth]{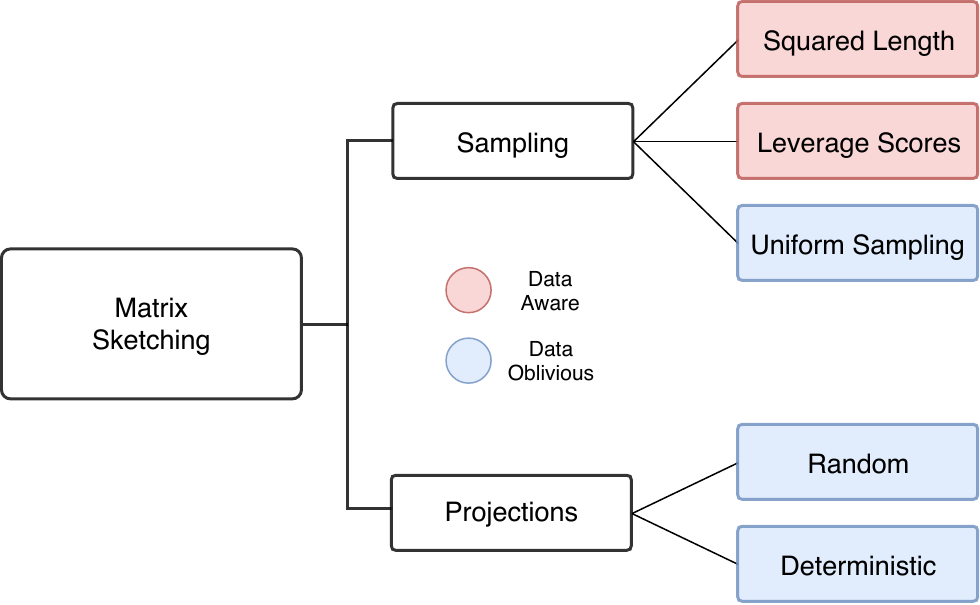}
	\vspace{-.05in}
	\caption{Sketching techniques for high-dimensional data; this work uses random projections.}
	\label{fig:sketching_techniques}
\end{figure}

Efforts to address the computational challenges of DMD can be traced back to the original paper by Schmid~\cite{Schmid2010jfm}. It was recognized that DMD could be computed in a lower dimensional space and then used to reconstruct the dominant modes and dynamics in the original high-dimensional space. This philosophy has been adopted in several strategies to determine the low-dimensional subspace, as shown in Fig.~\ref{fig:sketching_techniques}.  
Since then, efficient parallelized algorithms have been proposed for computations at scale~\cite{Sayadi2016tcfd}.  

In~\cite{Schmid2010jfm}, the high-dimensional data are projected onto the linear span of the associated POD modes, enabling an inexpensive low-dimensional DMD.  
The high-dimensional DMD modes have the same temporal dynamics as the low-dimensional modes and are obtained via lifting with the POD modes.
However, computing the POD modes may be expensive and strategies have been developed to alleviate this cost. An algorithm utilizing the randomized singular value decomposition (rSVD) was presented in~\cite{erichson2015} and~\cite{rdmdBistrian,bistrian2018efficiency}. While this approach is reliable and robust to noise, only the computation of the SVD is accelerated, so that subsequent computational steps involved in the DMD algorithm remain expensive. The paper by Bistrian and Navon~\cite{rdmdBistrian} has other advantages, including identifying the most influential DMD modes and guarantees that the low-order model satisfies the boundary conditions of the full model. 

As an alternative to using POD modes, a low-dimensional approximation of the data can be obtained by random projections~\cite{Brunton2015jcd}. 
This is justified by the Johnson-Lindenstrauss lemma~\cite{JL:1984,Fowler:2009} which provides statistical guarantees on the preservation of the distances between the data when projected into the low-dimensional space. This approach avoids the cost of computing POD modes and has favorable statistical error bounds. 

In contrast with these projection-based approaches, alternative strategies for reducing the dimension of the data involve computations with subsampled snapshots. In this case, the reduction of the computational cost comes from the fact that only a subset of the data is used to derive the DMD and the full data is never processed by the algorithm.  
DMD is then determined from a \emph{sketch-sampled} database. 
Variants of this approach differ in the sampling strategy; for example,~\cite{Brunton2015jcd} and~\cite{Erichson2016cdmd} develop the compressed dynamic mode decomposition, which involves forming a small data matrix by randomly projecting the high-dimensional row-space of a larger data matrix.  
This approach has been successful in computational fluid dynamics and video processing applications, where the data have relatively low noise.  
However, this is a suboptimal strategy, which leads to a large variance in the performance due to poor statistical guarantees in the resulting sketched data matrix. 
Clustering techniques can be employed to select a relevant subset of the data~\cite{Gueniat2015pof}. Further alternatives derive from algebraic considerations, such as the maximization of the volume of the submatrix extracted from the data (e.g., Q-DEIM \cite{QDEIM_Drmac,Manohar2017csm}) or rely on energy-based arguments such as leverage-score and length-squared samplings~\cite{frieze2004fast,Mahoney2011,RandNLA}.

\subsection{Assumptions, Limitations, Applications, and Extensions}\label{sec:limitations}
Within a short time, DMD has become a workhorse algorithm for extracting dominant low-dimensional oscillating patterns from high-dimensional data.  
One key benefit of DMD is that it is equally valid for experimental and numerical data, as it does not rely on knowledge of the governing equations.  
Although it may not be stated explicitly, it is often assumed that the data are periodic or quasi-periodic in nature.  
DMD typically fails for data that is non-stationary or that exhibits broadband or intermittent phenomena.  
DMD is also quite sensitive to noise~\cite{duke2012error,Dawson2016ef}, although several algorithms exist to de-bias DMD results in the presence of noisy data~\cite{Dawson2016ef,Hemati2017tcfd}. 
Although not a limitation, most applications of DMD involve high-dimensional systems that evolve on a low-dimensional attractor, in which case the singular value decomposition is used to determine the subspace.
Despite the assumptions and limitations above, DMD has been widely applied on a variety of systems in fluid mechanics~\cite{Rowley2009jfm,schmid2011applications,lusseyran2011flow,seena2011dynamic,basley2013space}, epidemiology~\cite{Proctor2015ih}, neuroscience~\cite{brunton2016extracting}, video modeling~\cite{Erichson2016jrtp,erichson2015,kutz2015multiVideo}, robotics~\cite{Berger2014ieee}, and plasma physics~\cite{Taylor2017arxiv}. 
Much of the success of DMD relates to its simple formulation in terms of linear algebra.  
The simplicity of DMD has enabled several extensions, including for control~\cite{Proctor2016siads}, multiresolution analysis~\cite{Kutz2016siads}, recursive orthogonalization of modes for Galerkin projection~\cite{Noack2016jfm}, the use of time-delay coordinates~\cite{Tu2014jcd,Arbabi2016arxiv,Das2017arxiv}, sparse identification of dynamic regimes~\cite{Kramer2015arxiv}, Bayesian formulations~\cite{Takeishi2017JCAI}.
  
Many of the applications and extensions above fundamentally rely on the dynamics evolving on a low-dimensional subspace that is well characterized by a fast decaying singular value spectrum.  Although this is not a fundamental limitation of DMD, it is often a basic assumption, and we rely on this low-rank structure here.

\subsection{Contribution of the Present Work}
This work presents a randomized DMD algorithm that enables the accurate and efficient extraction of spatiotemporal coherent structures from high-dimensional data.  
The algorithm scales with the intrinsic rank of the dynamics, which are assumed to be low-dimensional, rather than the  ambient measurement dimension.  
We demonstrate the effectiveness of this algorithm on two data sets of increasing complexity, namely the canonical fluid flow past a circular cylinder and the high-dimensional sea-surface temperature dataset.  
Moreover, we show that is possible to control the accuracy of the algorithm via oversampling and power iterations.  
In order to promote reproducible research, our open-source code is available at \href{https://github.com/erichson/ristretto}\texttt{https://github.com/erichson/ristretto}.

The remainder of the paper is organized as follows. Section~\ref{sec:technical} presents notation and background concepts used throughout the paper.  
Section~\ref{sec:framework} outlines the probabilistic framework used to compute the randomized DMD, which is presented in Sec.~\ref{sec:rdmd}.
Section~\ref{sec:results} provides numerical results to demonstrate the performance of the randomized DMD algorithm. Final remarks and an outlook are given in Sec.~\ref{sec:conclusion}.

\section[Technical Preliminaries]{Technical Preliminaries}\label{sec:technical}
We now establish notation and provide a brief overview of the singular value decomposition (SVD) and the dynamic mode decomposition (DMD).

\subsection[Notation]{Notation}
Vectors in $\mathbb{R}^{n}$ and $\mathbb{C}^{n}$ are denoted as bold lowercase letters ${\mathbf{x}=[x_1,x_2,...,x_n]}$. 
Both real $\mathbb{R}^{n\times m}$ and complex $\mathbb{C}^{n\times m}$ matrices are denoted by bold capitals $\mathbf{X}$, and its entry at the $i$-th row and $j$-th column is denoted as $\mathbf{X}(i,j)$.
The Hermitian transpose of a matrix is denoted as $\mathbf{X}^*$. The spectral or operator norm of a matrix is defined as the largest singular value $\sigma_{\rm max}(\mathbf{X})$ of the matrix $\mathbf{X}$, i.e., the square root of the largest eigenvalue $\lambda_{\rm max}$ of the positive-semidefinite matrix $\mathbf{X}^*\mathbf{X}$:
\begin{equation*}
\| \mathbf{X} \|_2 = \sqrt{\lambda_{\rm max}(\mathbf{X}^*\mathbf{X})} = \sigma_{\rm max}(\mathbf{X}).
\end{equation*} 
The Frobenius norm of a matrix $\mathbf{X}$ is the positive square root of the sum of the absolute squares of its elements, which is equal to the positive square root of the trace of $\mathbf{X}^*\mathbf{X}$
\begin{equation*}
\| \mathbf{X} \|_F = \sqrt{ \sum_{i=1}^n\sum_{j=1}^m|\mathbf{X}(i,j)|^2} = \sqrt{\textrm{trace}(\mathbf{X}^*\mathbf{X})}.
\end{equation*} 
The relative reconstruction error is $\|\mathbf{X}-\mathbf{\widehat{X}}\|_F / \|\mathbf{X}\|_F$, where $\mathbf{\widehat{X}}$ is an approximation to the matrix $\mathbf{X}$.
The column space (range) of $\mathbf{X}$ is denoted $\text{col}(\mathbf{X})$ and the row space is $\text{row}(\mathbf{X})$. 

\subsection[The Singular Value Decomposition]{The Singular Value Decomposition}

Given an ${n\times m}$ matrix $\mathbf{X}$, the `economic' singular value decomposition (SVD) produces the factorization
\begin{equation}\label{eq:svd}
\mathbf{X} = \mathbf{U}\mathbf{\Sigma}\mathbf{V}^* = \sum_{i}^{r} \mathbf{u}_i \mathbf{\sigma}_i\mathbf{v}_i^*,
\end{equation} 
where $\mathbf{U} = [ \mathbf{u}_1, ... , \mathbf{u}_r ] \in \mathbb{R}^{n \times r}$ and $\mathbf{V} = [ \mathbf{v}_1, ... , \mathbf{v}_r ] \in \mathbb{R}^{m\times r}$ are orthonormal, $\mathbf{\Sigma} \in \mathbb{R}^{r \times r}$ is diagonal, and $r=\text{min}(m,n)$. The left singular vectors in $\mathbf{U}$ provide a basis for the column space of $\mathbf{X}$, and the right singular vectors in $\mathbf{V}$ form a basis for the row space of $\mathbf{X}$. 
$\mathbf{\Sigma}$ contains the corresponding nonnegative singular values $\sigma_{1} \geq ... \geq  \sigma_{r} \geq 0$. 
Often, only the $k$ dominant singular vectors and values are of interest, resulting in the low-rank SVD: 
\begin{equation*}
\label{svdLR}
\mathbf{X}_k = \mathbf{U}_k\mathbf{\Sigma}_k\mathbf{V}^*_k = [ \mathbf{u}_1, \ldots , \mathbf{u}_k ] 
\textrm{diag}(\sigma_1,\dots,\sigma_k) [ \mathbf{v}_1, \ldots , \mathbf{v}_k ]^* = \sum_{i}^{k} \mathbf{u}_i \mathbf{\sigma}_i\mathbf{v}_i^*.
\end{equation*}

\paragraph[The Moore-Penrose Pseudoinverse]{The Moore-Penrose Pseudoinverse} Given the singular value decomposition ${\bX} = {\bU} \mathbf{\Sigma} {\bV}^*$, the pseudoinverse $\bX^\dag$ is computed as 
\begin{equation}\label{eq:pinv}
{\bX}^\dag := {\bV} \mathbf{\Sigma}^\dag {\bU}^* = \sum_{i}^{r} \mathbf{v}_i \mathbf{\sigma}_i^\dag\mathbf{u}_i^*,
\end{equation}
where $\mathbf{\Sigma}^\dag$ is here understood as $$\mathbf{\Sigma}^\dag = \mathrm{diag}\left(\sigma_i^\dag\right), \qquad 
\sigma_i^\dag = \left\{\begin{tabular}{ll}$\sigma_i^{-1}$ & if~$\sigma_i > 0$,\\ $0$ & otherwise,\end{tabular}\right. \quad \forall \, i.$$
We use the Moore-Penrose pseudoinverse in the following to provide a least squares solution to a system of linear equations.

\subsection[Dynamic Mode Decomposition]{Dynamic Mode Decomposition}

DMD is a dimensionality reduction technique, originally introduced in the field of fluid dynamics~\cite{Schmid2010jfm,Rowley2009jfm,Kutz2016book}. 
The method extracts spatiotemporal coherent structures from an ordered time series of snapshots $\bx_0,\bx_1,...,\bx_m \in \mathbb{R}^{n}$, separated in time by a constant step $\Delta t$. Specifically, the aim is to find the eigenvectors and eigenvalues of the time-independent linear operator $\mathbf{A}: \mathbb{R}^{n}\rightarrow\mathbb{R}^{n}$ that best approximates the map from a given snapshot $\mathbf{x}_{j}$ to the subsequent snapshot $\mathbf{x}_{j+1}$ as
\begin{equation}\label{eq:dmdvec}
\mathbf{x}_{j+1} \approx \mathbf{A}\mathbf{x}_{j}. 
\end{equation}
%
%
Following~\cite{Tu2014jcd}, the deterministic DMD algorithm proceeds by first separating the snapshot sequence $\bx_0,\bx_1,...,\bx_m$ into two overlapping sets of data
%
\begin{align}\label{eq:Xsnap}
{\bX}_L \!=\! \begin{bmatrix}
\vline & \vline & & \vline \\
{\bx}_0 & {\bx}_1 & \cdots & {\bx}_{m-1}\\
\vline & \vline & & \vline
\end{bmatrix}, \hspace{0.1in} {\bX}_R \!=\! \begin{bmatrix}
\vline & \vline & & \vline \\
{\bx}_1 & {\bx}_2 & \cdots & {\bx}_m\\
\vline & \vline & & \vline
\end{bmatrix}.
\end{align}
%
${\bX}_L \in {\mathbb{R}}^{n\times m}$ and ${\bX}_R \in {\mathbb{R}}^{n\times m}$ are called the left and right snapshot sequences. Equation~\eqref{eq:dmdvec} can then be reformulated in matrix notation as
\begin{equation}
{\bX}_R \approx {\bA}{\bX}_L.
\end{equation}

Many variants of the DMD have been proposed since its introduction, as discussed in \cite{Kutz2016book}. Here, we discuss the \emph{exact DMD} formulation of Tu et al.~\cite{Tu2014jcd}.

\subsubsection{Non-Projected DMD}

In order to find an estimate for the linear map $\bA$, the following least-squares problem can be formulated
\begin{equation}\label{eq:Xls}
\widehat{\bA} =	\argmin_{\bA} \|{\bX}_R - {\bA}{\bX}_L\|_{F}^{2}.
\end{equation}
The estimator for the best-fit linear map is given in closed form as
\begin{equation}\label{eq:AdirectPI}
\widehat{\bA} := \bX_R\bX_L^\dagger,
\end{equation}
where $\bX^\dag$ denotes the pseudoinverse from Eq.~\eqref{eq:pinv}.
The DMD modes are the eigenvectors of $\widehat{\bA} \in \mathbb{R}^{n\times n}$. 

\subsubsection{Projected DMD}\label{sec:projectedDMD}


If the data are high-dimensional (i.e., $n$ is large), the linear map $\widehat{\bA}$ in Eq.~\eqref{eq:AdirectPI} may be intractable to evaluate and analyze directly. Instead, a rank-reduced approximation of the linear map is determined by projecting it onto a $k$-dimensional subspace, $k \le \mathrm{min}\left(n ,m\right)$.
We denote $\bX_L = \bU \bSigma \bV^*$ the singular value decomposition of $\bX_L$.
The projected map onto the class of rank-$k$ linear operators can be obtained by pre- and postmultiplying Eq.~\eqref{eq:AdirectPI} with $\bP_k \in \mathbb{R}^{n \times k}$ and $\bT_k \in \mathbb{R}^{n \times k}$:
\begin{equation}\label{eq:projectedAtilde1}
	\widetilde{\bA} := \bP_k^* \widehat{\bA} \bT_k =  \bP_k^* \bX_R \bX_L^\dagger \bT_k = \bP_k^* \bX_R \bV \bSigma^\dag \bU^* \bT_k,
\end{equation}
where $\widetilde{\bf A} \in {\mathbb{R}}^{k\times k}$ is the projected map, i.e., a low-dimensional system  matrix.

The \emph{projected DMD}~\cite{Schmid2010jfm} relies on the hypothesis that the columns of ${\bA}{\bX}_L$ are in the linear span of ${\bX}_L$. Within this approximation, projection matrices $\bP$ and $\bT$ are chosen such that they span a subspace of the column space of ${\bX}_L$. A convenient choice is given by the dominant (i.e., associated with the largest singular values) left singular vectors $\bU_k$ of $\bX_L$ (i.e., POD modes) so that $\bP_k = \bT_k = \bU_k$, $\bU_k^* \bU_k = \bI_k$, with $k \le \mathrm{rank}\left(\bX_L\right)$, and
\begin{equation}\label{eq:projectedAtilde}
\widetilde{\bA} = \bU_k^* \bX_R \bV_k \bSigma_k^{-1},
\end{equation}
since ${\bU}^* {\bU}_k = \begin{bmatrix} \mathbf{I}_k \\ \boldsymbol{0} \end{bmatrix}$.

Note that when we use the POD modes ${\bU}$, one can interpret $\widehat{\bA} \bU_k$ as the shifted spatial POD modes over one time step $\Delta t$. Now, by premultiplying $\widehat{\bf A} \bU_k$ with ${\bU}_k^*$, one obtains a projected map $\widetilde{\bf A} = {\bU}_k^* \widehat{\bf A} \bU_k$ which encodes the cross-correlation of the POD modes ${\bU}_k$ with the time-shifted spatial POD modes. It thus becomes obvious that the DMD encodes more information about the temporal evolution of the system under consideration than the time-averaged POD modes~\cite{Schmid2010jfm}.

\begin{figure}[!b]
	\centering
		\begin{subfigure}[t]{0.31\textwidth}
		\centering
		\DeclareGraphicsExtensions{.pdf}
		\includegraphics[width=1\textwidth]{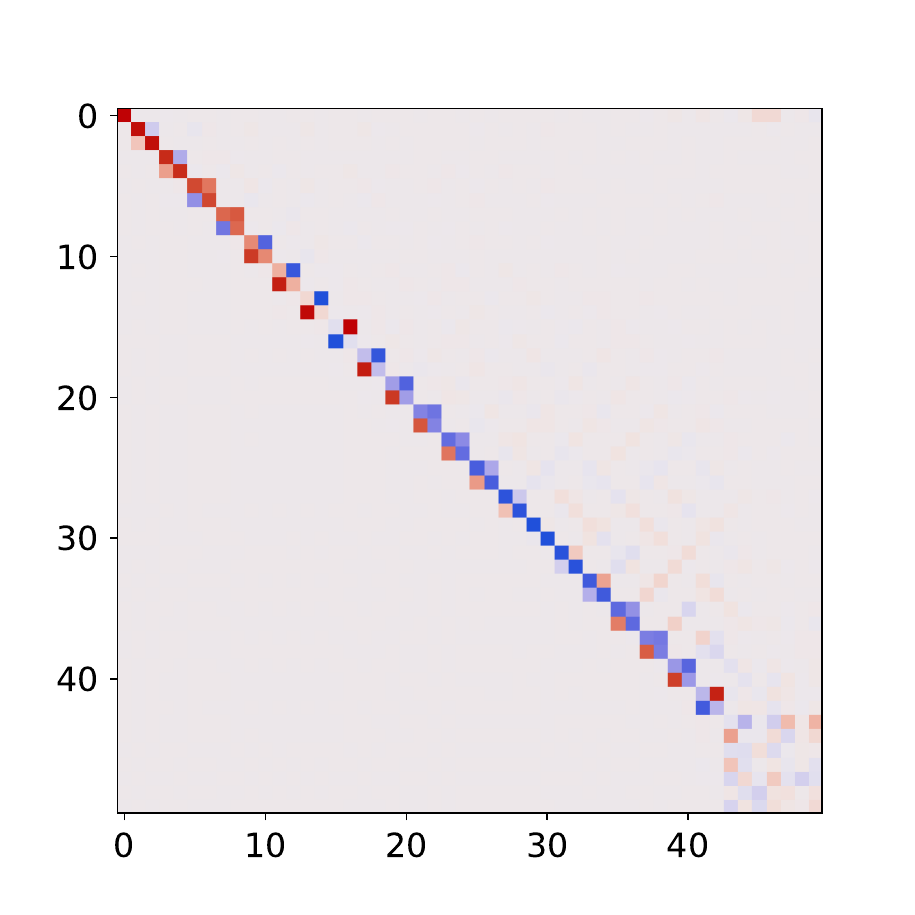}
		\caption{Non-projected linear map in the absence of white noise.}
	\end{subfigure}
	~	
	\begin{subfigure}[t]{0.31\textwidth}
		\centering
		\DeclareGraphicsExtensions{.pdf}
		\includegraphics[width=1\textwidth]{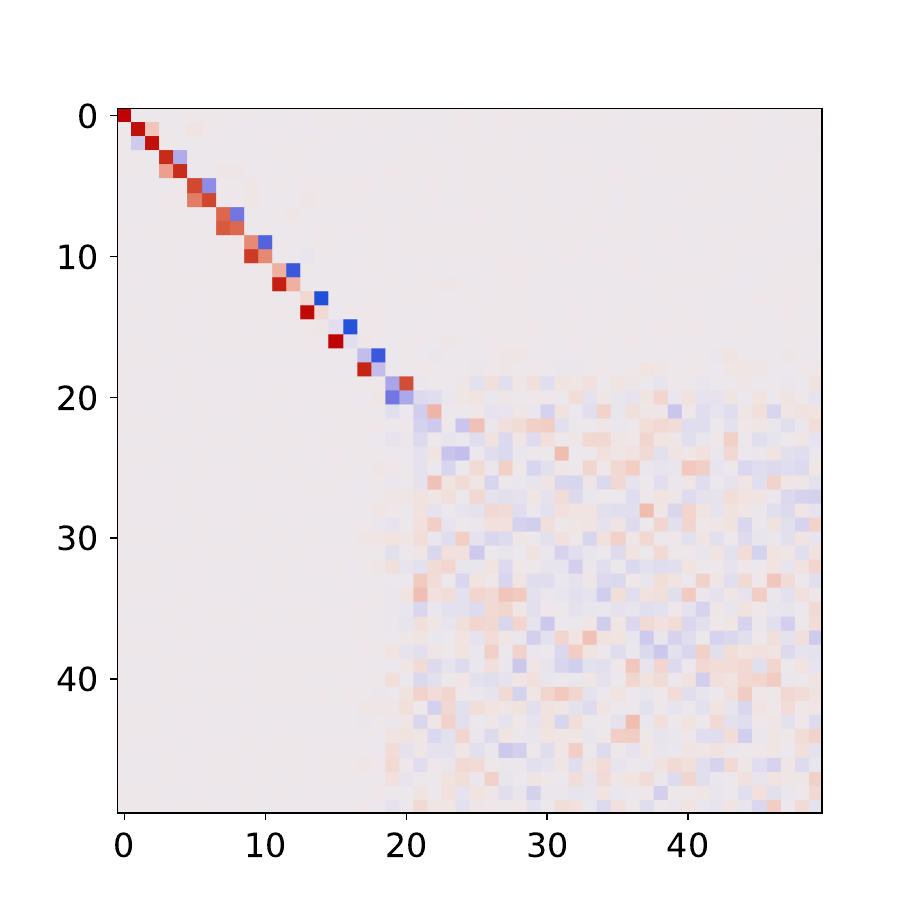}
		\caption{Non-projected linear map in the presence of white noise.}
	\end{subfigure}
	~
	\begin{subfigure}[t]{0.31\textwidth}
		\centering
		\DeclareGraphicsExtensions{.pdf}
		\includegraphics[width=1\textwidth]{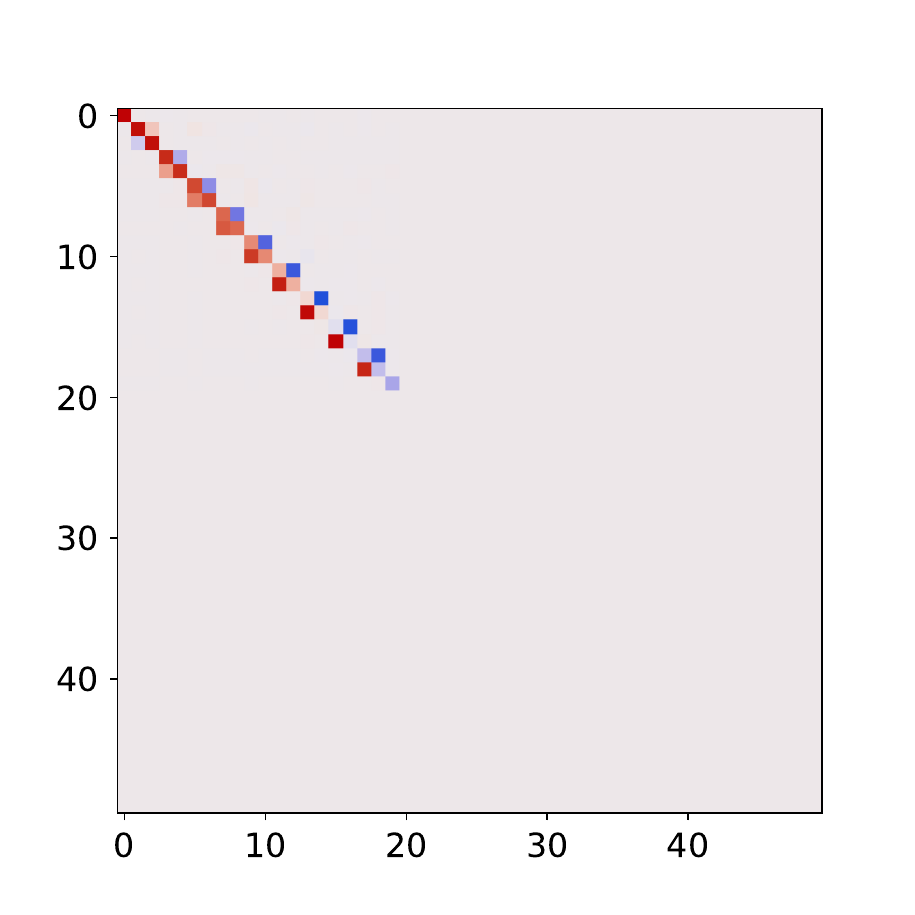}				
		\caption{Projected linear map in the presence of white noise.}
	\end{subfigure}		
	
	\caption{Illustration of the non-projected and projected map for toy data in the absence and presence of white noise. The regularization effect of the projected map reduces the influence of noise, while preserving the dominant information. }
	\label{fig:map}
\end{figure}

Low-dimensional projection is a form of spectral filtering which has the positive effect of dampening the influence of  noise~\cite{hansen2006deblurring,gorodnitsky1994analysis}. This effect is illustrated in Figure~\ref{fig:map}, which shows the non-projected linear map in absence and presence of white noise. The projected linear map avoids the amplification of noise and acts as a hard-threshold regularizer.
The difficulty is to choose the rank providing a good trade-off between suppressing the noise and retaining the useful information.
In practice, the optimal hard-threshold~\cite{gavish2014optimal} provides a good heuristic to determine the rank $k$. 

Once the linear map $\widetilde{\bA}$ is approximated, its eigendecomposition is computed	
\begin{equation}
	\widetilde{\bA}\widetilde{\bW} = \widetilde{\bW}\boldsymbol{\Lambda},  \label{eq:atildeweqwlamb}
\end{equation}
where the columns of $\widetilde{\bW} \in {\mathbb{C}}^{k\times k}$ are eigenvectors of $\widetilde{\bA}$, and $\boldsymbol{\Lambda} \in {\mathbb{C}}^{k\times k}$ is a diagonal matrix containing the corresponding eigenvalues $\lambda_j$. 
%
%
Eigenvalues of ${\bf \widetilde{\bA}}$ are eigenvalues of ${\bf \widehat{\bA}}$ and one may recover DMD modes as the columns of ${\bf W}_k$, \cite{Tu2014jcd}:
%
\begin{equation}
 {\bf W}_k := \bX_R \bV_k \bSigma_k^{-1} \widetilde{\bW}.
\end{equation}

In both the projected and the non-projected formulation, a singular value decomposition of $\bX_L$ is involved. 
This is computationally demanding, and the resources required to compute DMD can be tremendous for large data matrices.


\subsection{A Note on Regularization}

The projected algorithm described above relies on the truncated singular value decomposition (TSVD), also referred to as pseudo-inverse filter, to solve the unconstrained least-squares problem in Eq.~\eqref{eq:Xls}. 
Indeed, Figure~\ref{fig:map} illustrates that the low-rank approximation introduces an effective regularization effect. 
This warrants an extended discussion on regularization, following the work by Hansen~\cite{hansen1987truncatedsvd,hansen1990truncated,hansen2006deblurring,hansen1992modified}. 

In DMD, we often experience a linear system of equations $\bX_R = \bA \bX_L$ where $\bX_L$ or $\bX_R$ is ill-conditioned, i.e., the ratio between the largest and smallest singular value is large, so that the pseudo-inverse may magnify small errors. 
In this situation, regularization becomes crucial for computing an accurate estimate of $\bA$, since a small perturbation in $\bX_R$ or $\bX_L$ may result in a large perturbation in the solution. 
Tikhonov-Phillips regularization~\cite{tihonov1963solution,phillips1962technique}, also known as ridge regression in the statistical literature, is one of the most popular regularization techniques. The regularized least squares problem becomes
\begin{equation}
\widehat{\bA} =	\argmin_{\bA} \left\|{\bX}_R - {\bA}{\bX}_L\|_{F}^{2} + \lambda^2 \|\bA \right\|^2_F.
\end{equation}
More concisely, this can be expressed as the following augmented least-squares problem
\begin{equation}
\widehat{\bA} =	\argmin_{\bA} \left\| 	\begin{bmatrix}
{\bX}_R  \\
\mathbf{0}
\end{bmatrix} - {\bA} \begin{bmatrix}
{\bX}_L  \\
\lambda \mathbf{I} \end{bmatrix}  \right\|_{F}^{2}.
\end{equation}
The estimator for the map $\bA$ takes advantage of the regularized inverse 
\begin{equation}
\bX_\lambda^{\dagger} := [\bX^* \bX + \lambda^2 \mathbf{I}]^{-1} \bX^*= \bV \bSigma^+_\lambda \bU^* \quad \text{with} \quad \bSigma^+_\lambda := \diag\left(\frac{\sigma_1}{\sigma_1^2 + \lambda^2},...,\frac{\sigma_r}{\sigma_r^2 + \lambda^2} \right)
\end{equation}
where the additional regularization (controlled via the tuning parameter $\lambda$) improves the conditioning of the problem.
This form of regularization can also be seen as a smooth filter that attenuates the parts of the solution corresponding to the small singular values. Specifically, the filter $f$ takes the form~\cite{hansen1987truncatedsvd} 
\begin{equation}
f_i = \frac{\sigma_i^2}{\sigma_i^2 + \lambda^2}, \qquad i = 1, 2, ..., r.
\end{equation}
The Tikhonov-Phillips regularization scheme is closely related to the TSVD and the Wiener filter.
Indeed, the TSVD is known to be a method for regularizing ill-posed linear least-squares problems, which often produces very similar results~\cite{varah1979practical}.
More concretely, regularization via the TSVD can be seen as a hard-threshold filter, which takes the form
\begin{equation}
f_i = \begin{cases}
1, & \sigma_i \geq \sigma_k \\
0, & \sigma_i < \sigma_k 
\end{cases}.
\end{equation}
Here, $k$ controls how much smoothing (or low-pass filtering) is introduced, i.e., increasing $k$ includes more terms in the SVD expansion; thus components with higher frequencies are included. See~\cite{hansen2006deblurring} for further details. 
%
As a consequence of the regularization effect introduced by TSVD, the DMD algorithm requires a careful choice of the target-rank $k$ to compute a meaningful low-rank DMD approximation. Indeed, the solutions (i.e., the computed DMD modes and eigenvalues) may be sensitive to the amount of regularization which is controlled via $k$. In practice, one can treat the parameter $k$ as a tuning-parameter and find the optimal value via cross-validation.
\section[Probabilistic Framework]{Randomized Methods for Linear Algebra}\label{sec:framework}

%
 

With rapidly increasing volumes of measurement data from simulations and experiments, deterministic modal extraction algorithms may become prohibitively expensive.  
Fortunately, a wide range of applications produce data which feature low-rank structure, i.e., the rank of the data matrix, given by the number of independent columns or rows, is much smaller than the ambient dimension of the data space. 
In other words, the data feature a large amount of redundant information. 
In this case, randomized methods allow one to efficiently produce familiar decompositions from linear algebra. 
These so-called \emph{randomized} numerical methods have the potential to transform computational linear algebra, providing accurate matrix decompositions at a fraction of the cost of deterministic methods.
Indeed, over the past two decades, probabilistic algorithms have been prominent for computing low-rank matrix approximations, as described in a number of excellent surveys~\cite{halko2011rand,Mahoney2011,RandNLA,Brunton2018book}. 
The idea is to use randomness as a computational strategy to find a smaller representation, often denoted as \emph{sketch}. 
This sketch can be used to compute an approximate low-rank factorization for the high-dimensional data matrix ${\bX}\in \mathbb{R}^{n \times m}$. 

Broadly, these techniques can be divided into random sampling and random projection-based approaches. The former class of techniques carefully samples and rescales some `interesting' rows or columns to form a sketch. For instance, we can sample rows by relying on energy-based arguments such as leverage-scores and length-squared samplings~\cite{frieze2004fast,Mahoney2011,RandNLA}. The second class of random projection techniques relies on favorable properties of random matrices and forms the sketch as a randomly weighted linear combination of the columns or rows. Computationally efficient strategies to form such a sketch include the subsampled randomized Hadamard transform~\cite{sarlos2006improved,tropp2011improved} and the CountSketch~\cite{woodruff2014sketching}. Choosing between the different sampling and random projection strategies depends on the particular application. In general, sampling is more computational efficient, while random projections preserve more of the information in the data. 

\subsection{Probabilistic Framework} One of the most effective off-the-shelf methods to form a sketch is the probabilistic framework proposed in the seminal work by Halko et al.~\cite{halko2011rand}: 
\begin{itemize}
	\item \textbf{Stage A:} Given a desired target-rank $k$, find a near-optimal orthonormal basis $\mathbf{Q} \in \mathbb{R}^{n\times k} $ for the range of the input matrix $\bX$.
	\item \textbf{Stage B:} Given the near-optimal basis $\mathbf{Q}$, project the input matrix onto the low-dimensional space, resulting in $\mathbf{B} \in \mathbb{R}^{k\times m}$. This smaller matrix can then be used to compute a near-optimal low-rank approximation.
\end{itemize}

\subsubsection{Stage A: Computing a Near-Optimal Basis}
The first stage is used to approximate the range of the input matrix. Given a target rank $k \ll \min(m,n)$, the aim is to compute a near-optimal basis $\mathbf{Q} \in \mathbb{R}^{n\times k}$ for the input matrix $\mathbf{X} \in \mathbb{R}^{n\times m}$ such that
\begin{equation}\label{eq:AQQA}
\mathbf{X} \approx \mathbf{Q}\mathbf{Q}^*\mathbf{X}.
\end{equation}
Specifically, the range of the high-dimensional input matrix is sampled using the concept of random projections. Thus a basis is efficiently computed as
\begin{equation}\label{(eq:YXC)}
\mathbf{Y} = \mathbf{X}\mathbf{\Omega},
\end{equation}
where $\mathbf{\Omega} \in \mathbb{R}^{m\times k}$ denotes a random test matrix drawn from the normal Gaussian distribution. 
However, the cost of dense matrix multiplications can be prohibitive. 
The time complexity is order $\mathcal{O}(nmk)$.
To improve this scaling, more sophisticated random test matrices, such as the subsampled randomized Hadamard transform, have been proposed~\cite{sarlos2006improved,tropp2011improved}. 
Indeed, the time complexity can be reduced to $\mathcal{O}(nm\cdot \text{log}(k))$ by using a structured random test matrix to sample the range of the input matrix.

The orthonormal basis $\mathbf{Q} \in \mathbb{R}^{n\times k}$ is obtained via the QR-decomposition $\mathbf{Y}=\mathbf{Q}\mathbf{R}$.

\subsubsection{Stage B: Computing a Low-Dimensional Matrix}
We now derive a smaller matrix $\mathbf{B}$ from the high-dimensional input matrix $\mathbf{X}$. Specifically, given the near-optimal basis $\mathbf{Q}$, the matrix $\mathbf{X}$ is projected onto the low-dimensional space
\begin{equation}
\mathbf{B} = \mathbf{Q}^*\mathbf{X},
\end{equation}
which yields the smaller matrix $\mathbf{B} \in \mathbb{R}^{k\times m}$. This process preserves the geometric structure in a Euclidean sense, so that angles between vectors and their lengths are preserved.  It follows that 
\begin{equation}
\mathbf{X} \approx \mathbf{Q}\mathbf{B}.
\end{equation}

\subsubsection{Oversampling}
In theory, if the matrix $\mathbf{X}$ has exact rank $k$, the sampled matrix $\mathbf{Y}$ spans a basis for the column space, with high probability. In practice, however, it is common that the truncated singular values $\{\sigma_i\}_{i\geq k+1}$ are nonzero. Thus, it helps to construct a slightly larger test matrix in order to obtain an improved basis. Thus, we compute $\mathbf{Y}$ using an $\mathbf{\Omega} \in \mathbb{R}^{m\times l}$ test matrix instead, where $l=k+p$. Thus, the oversampling parameter $p$ denotes the number of additional samples. In most situations small values $p=\{5,10\}$ are sufficient to obtain a good basis that is comparable to the best possible basis~\cite{halko2011rand}. 
                                                  
\subsubsection{Power Iteration Scheme}\label{sec:power}
A second strategy to improve the performance is the concept of power iterations~\cite{rokhlin2009randomized,gu2015subspace}. In particular, a slowly decaying singular value spectrum of the input matrix can seriously affect the quality of the approximated basis matrix $\bQ$. Thus, the method of power iterations is used to preprocess the input matrix in order to promote a faster decaying spectrum. The sampling matrix $\mathbf{Y}$ is obtained as
\begin{equation}
\mathbf{Y} = \big( ( \mathbf{X}\mathbf{X}^* )^{q} \mathbf{X} \big)  \mathbf{\Omega}
\end{equation}
where $q$ is an integer specifying the number of power iterations. 
With $\mathbf{X} = \mathbf{U} \mathbf{\Sigma} \mathbf{V}^*$, one has $\mathbf{X}^{(q)}:=( \mathbf{X}\mathbf{X}^* )^{q} \mathbf{X} = \mathbf{U} \mathbf{\Sigma}^{2q + 1} \mathbf{V}^*$. 
Hence, for $q > 0$, the preprocessed matrix $\mathbf{X}^{(q)}$ has a relatively fast decay of singular values compared to the input matrix $\mathbf{X}$.
The drawback of this method is that additional passes over the input matrix are required. However, as few as $q=\{1,2\}$ power iterations can considerably improve the approximation quality, even when the singular values of the input matrix decay slowly.    

Algorithm~\ref{Alg:randomizedQB} describes this overall procedure for the randomized QB decomposition.

\begin{algorithm*}[t]
	\scalebox{.90 }{		
		\begin{minipage}{210mm}
			\begin{tabbing}
				\hspace{10mm} \= \hspace{5mm} \= \hspace{5mm} \= \hspace{45mm} \=\kill
				
				(1)  \> \> $l = k + p$ \> \> {\color{black}$\textrm{Slight oversampling.}$}\\[1mm]			
				
				(2)  \> \> $\mathbf{\Omega} = \texttt{rand}(m,l)$ \> \>  {\color{black}$\textrm{Generate random test matrix.}$}\\[1mm]
				
				(3)  \> \> $\mathbf{Y} = \bX \mathbf{\Omega}$ \> \> {\color{black}$\textrm{Compute sampling matrix.}$}\\[1mm]
				
				(4)  \> \> \textbf{for} $j = 1,\dots,q$ \> \> {\color{black}$\textrm{Power iterations (optional).}$} \\[1mm]
				
				(5)  \> \> \> $[\bQ, \sim] = \texttt{qr}(\mathbf{Y})$ \> \\[1mm] 
				
				(6)  \> \> \> $[\mathbf{Z}, \sim] = \texttt{qr}(\bX^{*} \bQ)$ \\[1mm]
				
				(7)  \> \> \> $\mathbf{Y} = \bX \mathbf{Z}$ \\[1mm]
				
				(8)  \> \> \textbf{end for}\\[1mm]
				
				(9)  \> \> $[\bQ, \sim] = \texttt{qr}(\mathbf{Y})$ \> \> {\color{black}\textrm{Orthonormalize sampling matrix.}}
				\\[1mm]
				
				(10)  \> \> $\bB =  \bQ^{*} \bX$ \> \> {\color{black}$\textrm{Project input matrix to smaller space.}$} 						
			\end{tabbing}
	\end{minipage}}
	\centering
	\caption{Randomized QB decomposition.\\ Given a matrix $\mathbf{X} \in \mathbb{R}^{n \times m}$ and a target rank $k \ll \min(m,n)$, the basis matrix $\mathbf{Q} \in {\mathbb{R}}^{n\times k}$ with orthonormal columns and the smaller matrix $\mathbf{B} \in {\mathbb{R}}^{k\times m}$ are computed. The approximation quality can be controlled via oversampling $p$ and the computation of $q$ power iterations.}
	\label{Alg:randomizedQB}
\end{algorithm*}
%

\begin{remark}
	A direct implementation of the power iteration scheme, as outlined in Section~\ref{sec:power}, is numerically unstable due to round-off errors. Instead, the sampling matrix $\bY$ is orthogonalized between each computational step to improve the stability~\cite{halko2011rand}. 	
\end{remark}	

\begin{remark} As default values for the oversampling and power iteration parameter, we suggest $p=10$, and $q=2$. 
\end{remark}

\subsubsection{Theoretical Performance}
Both the concept of oversampling and the power iteration scheme provide control over the quality of the low-rank approximation. The average-case error behavior is described in the probabilistic framework as~\cite{martinsson2016randomized2}:
\begin{equation*}
\E_{\mu_{\bQ}} \| \bX - \bQ\bQ^*\bX  \|_F  \leq  \Bigg[ 1 + \sqrt{\frac{k}{p-1}} + \frac{e\sqrt{l}}{p} \cdot \sqrt{\text{min}(m,n)-k}\Bigg]^\frac{1}{2q+1} \sigma_{k+1}(\bX),
\end{equation*}
with $e \equiv \exp\left(1\right)$, $\E_{\mu_{\bQ}}$ the expectation operator over the probability measure of $\bQ$ and $\sigma_{k+1}(\bX)$ the $\left(k+1\right)$-th element of the decreasing sequence of singular values of $\bX$.
Here it is assumed that $p\geq2$. 
Thus, both oversampling and the computation of additional power iterations drive the approximation error down.

\subsubsection{Computational Considerations}
The steps involved in computing the approximate basis $\mathbf{Q}$ and the low-dimensional matrix $\mathbf{B}$ are simple to implement, and embarrassingly parallelizable. Thus, randomized algorithms can benefit from modern computational architectures, and they are particularly suitable for GPU-accelerated computing. Another favorable computational aspect is that only two passes over the input matrix are required in order to obtain the low-dimensional matrix. Pass efficiency is a crucial aspect when dealing with massive data matrices which are too large to fit into fast memory, since reading data from the hard disk is prohibitively slow and often constitutes the actual bottleneck.    
%

\subsection{Blocked Randomized Algorithm}
When dealing with massive fluid flows that are too large to read into fast memory, the extension to sequential, distributed, and parallel computing might be inevitable~\cite{Sayadi2016tcfd}. In particular, it might be necessary to distribute the data across processors which have no access to a shared memory to exchange information.
To address this limitation, Martinsson and Voronin~\cite{martinsson2016randomized} proposed a blocked scheme to compute the QB decomposition on smaller blocks of the data. 
The basic idea is that a given high-dimensional sequence of snapshots $\bx_0,\bx_1,...,\bx_m \in \mathbb{R}^{n}$ is subdivided into $b$ smaller blocks along the rows. The submatrices can then be processed in $b$ independent streams of calculations. Here, $b$ is assumed to be a power of two, and zero padding can be used in order to divide the data into blocks of the same size.
This scheme constructs the smaller matrix $\bB$ in a hierarchical fashion, for more details see \cite{voronin2015rsvdpack}.

In the following, we describe a modified and simple scheme that is well suited for practical applications such as computing the dynamic mode decomposition for large-scale fluid flows. 
Unlike~\cite{martinsson2016randomized,voronin2015rsvdpack}, which discuss a fixed-precision approximation, we focus on a fixed-rank approximation problem.
Further, our motivation is that it is often unnecessary to use the full hierarchical scheme if the desired target rank $k$ is relatively small. By small we mean that we can fit a matrix of dimension $(b\cdot k) \times m$ into fast memory. To be more precise, suppose again that we are given an input matrix $\mathbf{X}$ with $n$ rows and $m$ columns. 
We partition $\mathbf{X}$ along the rows into $b$ blocks $\mathbf{X}_i$ of dimension $n/b \times m$. To illustrate the scheme we set $b=4$ so that
\begin{equation}
{\bX} =
	\begin{bmatrix}
	{\bX}_1  \\
	{\bX}_2  \\
	{\bX}_3  \\
	{\bX}_4 
	\end{bmatrix}.
\end{equation} 
Next, we approximate each block $\mathbf{X}_i$ by a fix rank-$k$ approximation, using the QB decomposition described in Algorithm~\ref{Alg:randomizedQB}, which yields 
%
%
%
%
\begin{equation}
{\bX} \approx
\begin{bmatrix}
{\bQ}_1 {\bB}_1  \\
{\bQ}_2 {\bB}_2  \\
{\bQ}_3 {\bB}_3  \\
{\bQ}_4 {\bB}_4     
\end{bmatrix} =
\diag\left({\bQ}_1, {\bQ}_2, {\bQ}_3, {\bQ}_4\right) 
\begin{bmatrix}
{\bB}_1 	\\
{\bB}_2 	\\
{\bB}_3 	\\
{\bB}_4    
\end{bmatrix}.
\end{equation} 
We can then collect all matrices ${\bB_i} \in \mathbb{R}^{k \times m}$ and stack them together as
\begin{equation}
{\bK} =
\begin{bmatrix}
{\bB}_1^*&  	
{\bB}_2^*&  	
{\bB}_3^*&  	
{\bB}_4^*     
\end{bmatrix}^{*}.
\end{equation} 
Subsequently, we compute the QB decomposition of $\bK \in \mathbb{R}^{(b\cdot k) \times m}$ and obtain 
\begin{equation}
{\bK} \approx \widehat{\bQ} \bB.
\end{equation} 
The small matrix $\bB \in \mathbb{R}^{k \times m}$ can then be used to compute the randomized dynamic mode decomposition as described below. The basis matrix ${\bQ} \in \mathbb{R}^{n \times k}$ can be formed as
%

\begin{equation}
{\bQ} =
\diag\left({\bQ}_1, {\bQ}_2, {\bQ}_3, {\bQ}_4\right) \widehat{\bQ}.
\end{equation}

In practice, we choose the target-rank $l$ for the approximation slightly larger than $k$, i.e., $l=k+p$.

\section[Randomized Dynamic Mode Decomposition]{Randomized Dynamic Mode Decomposition}\label{sec:rdmd}

In many cases, even with increased measurement resolution, the data may have dominant coherent structures that define a low-dimensional attractor~\cite{HLBR_turb}; in fact, the presence of these structures was the original motivation for methods such as DMD.  
%
%
Hence, it seems natural to use randomized methods for linear algebra to accelerate the computation of the approximate low-rank DMD.

\subsection{Problem Formulation Revisited} 

We start our discussion by formulating the following least-squares problem to find an estimate ${\bf \widehat{A}_\bB}$ for the projected linear map $\widetilde{\bf A} \in {\mathbb{R}}^{l\times l}$
in terms of the projected snapshot sequences ${\bB}_L := \bP^* {\bX}_L \in {\mathbb{R}}^{l\times m}$ and ${\bB}_R := \bP^* {\bX}_R \in {\mathbb{R}}^{l\times m}$:
\begin{equation}\label{eq:projectedLS}
{\bf \widehat{A}_\bB} := \argmin_{\bA_\bB} \|{\bB}_R - {\bA_\bB}{\bB_L}\|_{F}^{2},
\end{equation}
where $\bP \in {\mathbb{R}}^{n\times l}$ is a projection matrix discussed below. Recall that $l = k+p$, where $k$ denotes the desired target rank, and $p$ is the oversampling parameter.
Here we make the assumption that $\text{col}(\mathbf{\bP}) \approx \text{col}(\mathbf{X}_L)$ as well as $\text{col}(\mathbf{\bP}) \approx \text{col}(\mathbf{X}_R)$.
The question remains how to construct $\bP$ in order to quickly compute ${\bB}_L$ and ${\bB}_R$. The dominant left singular vectors of ${\bX}_R$ provide a good choice; however, the computational demands to compute the deterministic SVD may be prohibitive for large data sets. 
Next, we discuss randomized methods to compute DMD.

\subsection{Compressed Dynamic Mode Decomposition}

Brunton et al.~\cite{Brunton2015jcd} proposed one of the first algorithms using the idea of matrix sketching to find small matrices ${\bB}_L$ and ${\bB}_R$.
The algorithm proceeds by forming a small number of random linear combinations of the rows of the left and right snapshot matrices to form the representation
\begin{equation}\label{(eq:YCX)}
\mathbf{B}_L = \mathbf{S}\mathbf{X}_L, \quad\mathbf{B}_R = \mathbf{S}\mathbf{X}_R.
\end{equation}
with $\mathbf{S}\in\mathbb{R}^{l\times n}$ a random test matrix. As discussed in Section~\ref{sec:framework}, $\mathbf{S}$ can be constructed by drawing its entries from the standard normal distribution. While using a Gaussian random test matrix has beneficial theoretical properties, the dense matrix multiplication $\mathbf{S}\mathbf{X}_L$ can be expensive for large dense matrices. 
Alternatively, $\mathbf{S}$ can be chosen to be \textit{a random and rescaled subset} of the identity matrix, i.e., $\mathbf{B}$ is formed by sampling and rescaling $l$ rows from $\mathbf{X}$ with probability $p_i$.  
Thus, $\mathbf{S}$ is very sparse and is not required to be explicitly constructed and stored.
More concretely, we sample the $i$th row of $\mathbf{X}$ with probability $p_i$ and, if sampled, the row is rescaled by the factor $1/\sqrt{l \cdot p_i}$ in order to yield an unbiased estimator~\cite{drineas2017lectures}.
A naive approach is to sample rows with uniform probability $p_i = 1/n$. 
Uniform random sampling may work if the information is evenly distributed across the input data, so that dominant structures are not spatially localized.
Indeed, it was demonstrated in~\cite{Brunton2015jcd} that the dominant DMD modes and eigenvalues for some low-rank fluid flows can be obtained from a massively undersampled sketch.
However, the variance can be large and the performance may be poor in the presence of white noise. 
While not explored by~\cite{Brunton2015jcd}, the performance can be improved using more sophisticated sampling techniques, such as leverage-score sampling~\cite{drineas2012fast,ma2015statistical,drineas2017lectures}.
Better performance is expected when using structured random test matrices such as the subsampled randomized Hadamard transform~\cite{sarlos2006improved,tropp2011improved} or the CountSketch~\cite{woodruff2014sketching}. Both of these methods enable efficient matrix multiplications, yet they are more computationally demanding than random sampling.

The shortcoming of the algorithm in~\cite{Brunton2015jcd} is that $l$ depends on the ambient dimension of the measurement space of the input matrix. Thus, even if the data matrix is low rank, a large $l$ may be required to compute the approximate DMD. Next, we discuss a randomized algorithm which depends on the intrinsic rank of the input matrix.

\subsection{An Improved Randomized Scheme}\label{sec:rdmdpres}

In the following, we present a novel algorithm to compute the dynamic mode decomposition using the probabilistic framework above. The proposed algorithm is simple to implement and can fully benefit from modern computational architectures, e.g., parallelization and multithreading.

Given a sequence of snapshots $\bx_0,\bx_1,...,\bx_m \in \mathbb{R}^{n}$, we first compute the near-optimal basis $\mathbf{Q} \in \mathbb{R}^{n\times l}$ using the randomized methods which we have discussed in Section~\ref{sec:framework}. Then, we project the data onto the low-dimensional space, so that we obtain the low-dimensional sequence of snapshots 
\[
\bb_0,\bb_1,...,\bb_m := \mathbf{Q}^{*} \bx_0,  \mathbf{Q}^{*} \bx_1,...,  \mathbf{Q}^{*}\bx_m \in \mathbb{R}^{l}. 
\]
More concisely, we can express this as
\[
\bB :=  \mathbf{Q}^{*} \bX. 
\]
Next, we separate the sequence $\bb_0,\bb_1,...,\bb_m$ into two overlapping matrices ${\bB}_L \in {\mathbb{R}}^{l\times m}$ and ${\bB}_R \in {\mathbb{R}}^{l\times m}$
%
\begin{align}
{\bB}_L \!=\! \begin{bmatrix}
\vline & \vline & & \vline \\
{\bb}_0 & {\bb}_1 & \cdots & {\bb}_{m-1}\\
\vline & \vline & & \vline
\end{bmatrix}, \hspace{0.1in} {\bB}_R \!=\! \begin{bmatrix}
\vline & \vline & & \vline \\
{\bb}_1 & {\bb}_2 & \cdots & {\bb}_m\\
\vline & \vline & & \vline
\end{bmatrix}.
\end{align}
This leads to the following projected least-squares problem
\begin{equation}
\widehat{\bA}_\bB =	\argmin_{\bA_\bB} \|{\bB}_R - {\bA_\bB}{\bB_L}\|_{F}^{2}.
\end{equation}
%
Using the pseudoinverse, the estimator for the linear map $\widehat{\bA}_\bB\in {\mathbb{R}}^{l\times l}$ is defined as
\begin{equation}\label{eq:defMb}
	\mathbf{\widehat{A}}_{\mathbf{B}} := \bB_R \bB_L^\dagger = \bB_R\mathbf{V}\bSig^{-1}\mathbf{\widetilde{U}}^*,
\end{equation}
%
where $\widetilde{\bU}\in {\mathbb{R}}^{l\times k}$ and ${\bV}\in {\mathbb{R}}^{m\times k}$ are the truncated left and right singular vectors of $\bB_L$.
%
The diagonal matrix $\bSig \in {\mathbb{R}}^{k\times k}$ contains the corresponding singular values.
If $p=0$, then $l=k$, and no truncation is needed.
Next, $\mathbf{\widehat{A}}_{\mathbf{B}}$ is  projected onto the left singular vectors
\begin{subequations}
	\begin{align}
	\widetilde{\mathbf{A}}_{\mathbf{B}} &= \widetilde{\mathbf{U}}^*\mathbf{\mathbf{\widehat{A}}_{\mathbf{B}}}\mathbf{\widetilde{U}}\\
	&= \mathbf{\widetilde{U}}^*\bB_R\mathbf{V}\bSig^{-1}.
	\end{align}
\end{subequations}
The DMD modes, containing the spatial information, are then obtained by computing the eigendecomposition of $\widetilde{\mathbf{A}}_{\mathbf{B}}\in {\mathbb{R}}^{k\times k}$
\begin{equation} \label{eq:eigenAtildeB}
\widetilde{\mathbf{A}}_{\bB}\widetilde{\mathbf{W}}_{\bB} = \widetilde{\mathbf{W}}_{\bB}\boldsymbol{\Lambda},
\end{equation}
where the columns of ${\widetilde{\bW}_\bB}\in {\mathbb{C}}^{k\times k}$ are eigenvectors $\widetilde{\mathbf{w}}_{{\bB},j}$, and $\boldsymbol{\Lambda}\in {\mathbb{C}}^{k\times k}$ is a diagonal matrix containing the corresponding eigenvalues $\lambda_j$.
%
The high-dimensional DMD modes $\mathbf{W} \in {\mathbb{C}}^{n\times k}$  may be recovered as
\begin{equation}\label{Eq:rDMDModes}
\mathbf{W} = \mathbf{Q}\bB_R\mathbf{V}\bSig^{-1}\widetilde{\bW}_\bB.
\end{equation}
%
%
%
The computational steps for a practical implementation are sketched in Algorithm~\ref{Alg:rDMDalgorithm}. 
\begin{algorithm*}[!b]
	\scalebox{.90 }{		
		\begin{minipage}{210mm}
			\begin{tabbing}
				\hspace{10mm} \= \hspace{5mm} \= \hspace{5mm} \= \hspace{45mm} \=\kill
						
				(1)  \> \> $[\bQ, \bB] = \texttt{rqb}(\mathbf{X}, p, q)$ \> \> {\color{black}\textrm{Randomized QB decomposition (Alg.~\ref{Alg:randomizedQB}).}}
				\\[1mm]
				
				(2)  \> \> $\bB \rightarrow \left\{\bB_L,\bB_R\right\} $ \> \> {\color{black}\textrm{Left/right low-dimensional snapshot matrix.}} \\[1mm]		
				
				(3)  \> \> $[\widetilde{\mathbf{U}},{\bSig},{\mathbf{V}}] = \texttt{svd}(\bB_L, k)$ \> \> {\color{black}\textrm{Truncated SVD.}}
				\\[1mm]
				
				(4)  \> \> $\widetilde{\mathbf{A}}_{\bB} = \widetilde{\mathbf{U}}^{*}  \bB_R  \mathbf{V}  \bSig^{-1}$ \> \> {\color{black}\textrm{Least squares fit.}} \\[1mm]
				
				(5)  \> \> $[\mathbf{\widetilde{W}}_{\bB},{\mathbf{\Lambda}}] = \texttt{eig}(\mathbf{\widetilde{\mathbf{A}}_{\bB}})$  \> \> {\color{black}\textrm{Eigenvalue decomposition.}} \\[1mm]
				
				(6)  \> \> {$\mathbf{W} \gets \mathbf{Q} \bB_R \mathbf{V} \bSig^{-1} \mathbf{\widetilde{W}}_{\bB}$}  \> \> {\color{black}\textrm{Recover high-dimensional DMD modes $\mathbf{W}$.}} 
				
				
			\end{tabbing}
	\end{minipage}}
	\centering
	\caption{Randomized Dynamic Mode Decomposition (rDMD).\\ Given a snapshot matrix $\mathbf{X}$ and a target rank $k$, the near-optimal dominant dynamic modes $\mathbf{W}$ and eigenvalues $\mathbf{\Lambda}$ are computed. The approximation quality can be controlled via oversampling $p$ and the computation of power iterations $q$. An implementation in {Python} is available via the GIT repository \href{https://github.com/erichson/ristretto}\texttt{https://github.com/erichson/ristretto}.}.  
	\label{Alg:rDMDalgorithm}
\end{algorithm*}
\subsubsection{Derivation}
In the following we outline the derivation of Eq.~\eqref{Eq:rDMDModes}.
Recalling Eq.~\eqref{eq:AdirectPI}, the estimated transfer operator $\mathbf{\widehat{A}}$ is defined as
\begin{equation}\label{eq:Aexactapp}
\widehat{\bA} := \bX_R \bX_L^\dagger.
\end{equation}
%
Random projections of the data matrix $\mathbf{X}$ result in an approximate orthonormal basis $\bQ$ for the range of $\mathbf{\widehat{A}}$, as described in Sec.~\ref{sec:framework}. The sampling strategy is assumed to be efficient enough so that $\bX_L \approx \bQ \bQ^* \bX_L$ and $\bX_R \approx \bQ \bQ^* \bX_R$. Equation~\eqref{eq:Aexactapp} then becomes
\begin{equation}\label{eq:AexactappQ}
\widehat{\bA} \approx \left(\bQ \bQ^* \bX_R\right) \left(\bQ \bQ^* \bX_L\right)^\dagger.
\end{equation}
Letting $\bB_L:=\bQ^*\bX_L$ and $\bB_R:=\bQ^*\bX_R$, the projected transfer operator estimate is defined as in Eq.~\eqref{eq:defMb}, $\mathbf{\widehat{A}}_{\mathbf{B}} := \bB_R \bB_L^\dagger$. Substituting this into Eq.~\eqref{eq:AexactappQ} leads to
\begin{equation}\label{eq:eigenhatA}
\widehat{\bA} \approx \bQ \mathbf{\widehat{A}}_{\mathbf{B}} \bQ^*.
\end{equation}
Letting $\bB_L = \widetilde{\bU} \bSigma \bV^*$ be the SVD of $\bB_L$, and left- and right-projecting $\widehat{\bA}_{\bB}$ onto the dominant $k$ left singular vectors $\widetilde{\bU}$, yields
\begin{equation}\label{eq:eigentildeA}
\widetilde{\bA}_{\bB} := \widetilde{\bU}^* \mathbf{\widehat{A}}_{\mathbf{B}} \widetilde{\bU}.
\end{equation}
The eigendecomposition is given by $\widetilde{\bA}_{\bB} \widetilde{\bW}_{\bB} = \widetilde{\bW}_{\bB} {\boldsymbol{\Lambda}}$, as in Eq.~\eqref{eq:eigenAtildeB}.

Let $\widehat{\bW}_{\bB} := \bB_R \, \bV \bSigma^{-1} \widetilde{\bW}_{\bB}$. It is simple to verify that this is an eigenvector of $\widehat{\bA}_{\bB}$:
\begin{eqnarray}
\widehat{\bA}_{\bB} \widehat{\bW}_{\bB} & = & \bB_R \, \bB_L^\dagger \bB_R \, \bV \bSigma^{-1} \widetilde{\bW}_{\bB}, \nonumber \\
& = & \bB_R \, \bV \bSigma^{-1} \widetilde{\bA}_{\bB} \widehat{\bW}_{\bB}, \nonumber \\
& = & \bB_R \, \bV \bSigma^{-1} \widetilde{\bW}_{\bB} \boldsymbol{\Lambda}, \nonumber \\
& = & \widehat{\bW}_{\bB} \boldsymbol{\Lambda}. \nonumber
\end{eqnarray}
Substituting the eigendecomposition of $\widehat{\bA}_{\bB}$ in Eq.~\eqref{eq:eigenhatA} and right-multiplying by $\bQ \widehat{\bW}_{\bB}$ leads to
\begin{equation}
\widehat{\bA} \bQ \widehat{\bW}_{\bB} \approx \bQ \widehat{\bW}_{\bB} \boldsymbol{\Lambda},
\end{equation}
so that identification with $\widehat{\bA} \bW = \bW \widetilde{\boldsymbol{\Lambda}}$ verifies the claim in Eq.~\eqref{Eq:rDMDModes} on the eigendecomposition of $\widehat{\bA}$:
\begin{equation*}
{\bW} \approx \bQ \bB_R \, \bV \bSigma^{-1} \widetilde{\bW}_{\bB}, \qquad \widetilde{\boldsymbol{\Lambda}} \approx \boldsymbol{\Lambda}.
\end{equation*}

\section[Results]{Numerical Results}\label{sec:results}

In the following we present numerical results demonstrating the performance of the randomized dynamic mode decomposition (rDMD). First, we provide some visual results for a flow behind a cylinder and climate data. Then, we evaluate and compare the computational time and accuracy of randomized algorithm to the deterministic algorithm. 

All computations are performed using Amazon Web Services (AWS). We use a G3 instance with 16 Xeon E5-2686 CPUs and 1 NVIDIA Tesla M60 GPU.   The underlying numerical linear algebra routines are accelerated using the Intel Math Kernel Library (MKL).

\subsection{Fluid Flow Behind a Cylinder}
As a canonical example, we consider the fluid flow behind a cylinder at Reynolds number $\Rey=100$. The data consist of a sequence of $151$ snapshots of fluid vorticity fields on a $449\times 199$ grid\footnote{Data for this example may be downloaded at dmdbook.com/DATA.zip.}, computed using the immersed boundary projection method~\cite{taira2007ibpm,colonius2008ibpm}. The flow features a periodically shedding wake structure and the resulting dataset is low rank. While this data set poses no computational challenge, it demonstrates the accuracy and quality of the randomized approximation on an example that builds intuition.
Flattening and concatenating the snapshots horizontally yields a matrix of dimension $\mathbf{X} \in \mathbb{R}^{89,351 \times 151}$, i.e., the columns are the flattened snapshots $\mathbf{x} \in \mathbb{R}^{449 \times 199}$.

\begin{figure}[!t]
	\centering
	\begin{subfigure}[t]{0.45\textwidth}
		\centering
		\DeclareGraphicsExtensions{.pdf}
		\includegraphics[width=1\textwidth]{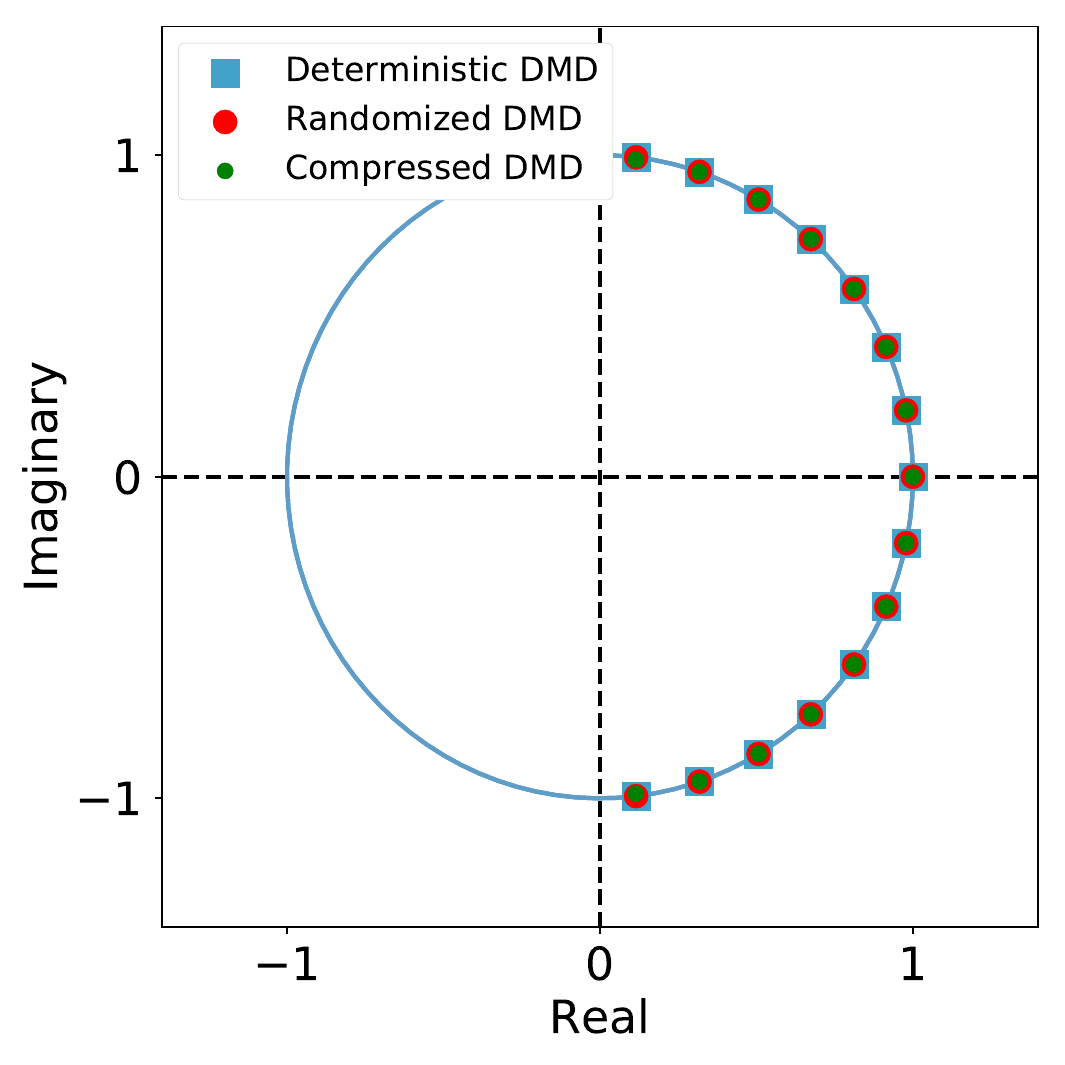}
		\caption{No noise.}
		\label{fig:flow_eigA}
	\end{subfigure}
	~
	\begin{subfigure}[t]{0.45\textwidth}
		\centering
		\DeclareGraphicsExtensions{.pdf}
		\includegraphics[width=1\textwidth]{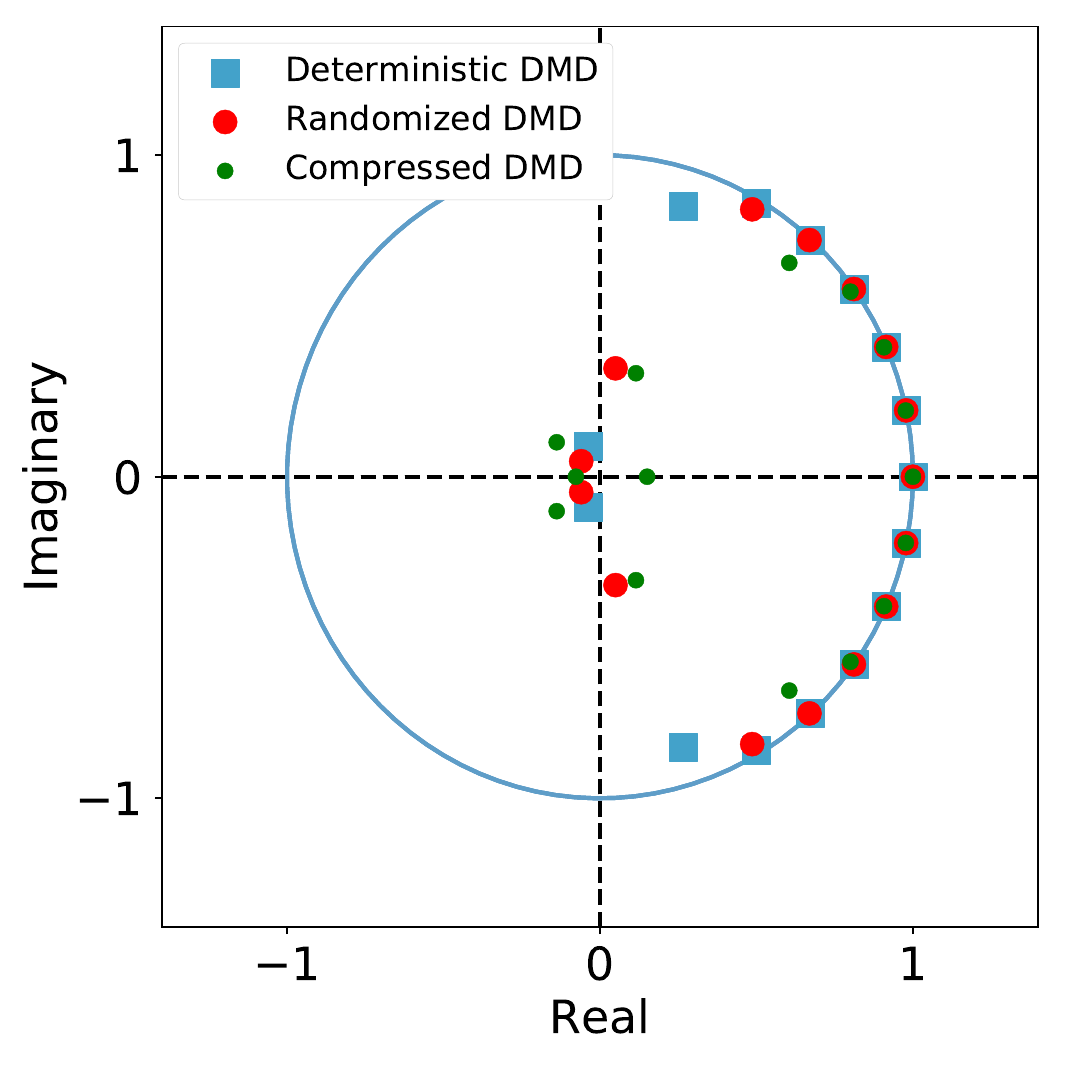}				
		\caption{Noisy with SNR of 10.}
		\label{fig:flow_eigB}
	\end{subfigure}		
	\vspace{-.2in}
	\caption{DMD eigenvalues for fluid flow behind a cylinder. Both the compressed and randomized DMD algorithms capture the eigenvalues in the absence of noise (a). In the presence of white noise with SNR of 10, rDMD performs better than sampling rows (b).}
	\label{fig:flow_eig}
\end{figure}

\begin{figure}[!t]
	\centering
	\begin{subfigure}[t]{0.48\textwidth}
		\centering
		\DeclareGraphicsExtensions{.png}
		\includegraphics[width=1\textwidth]{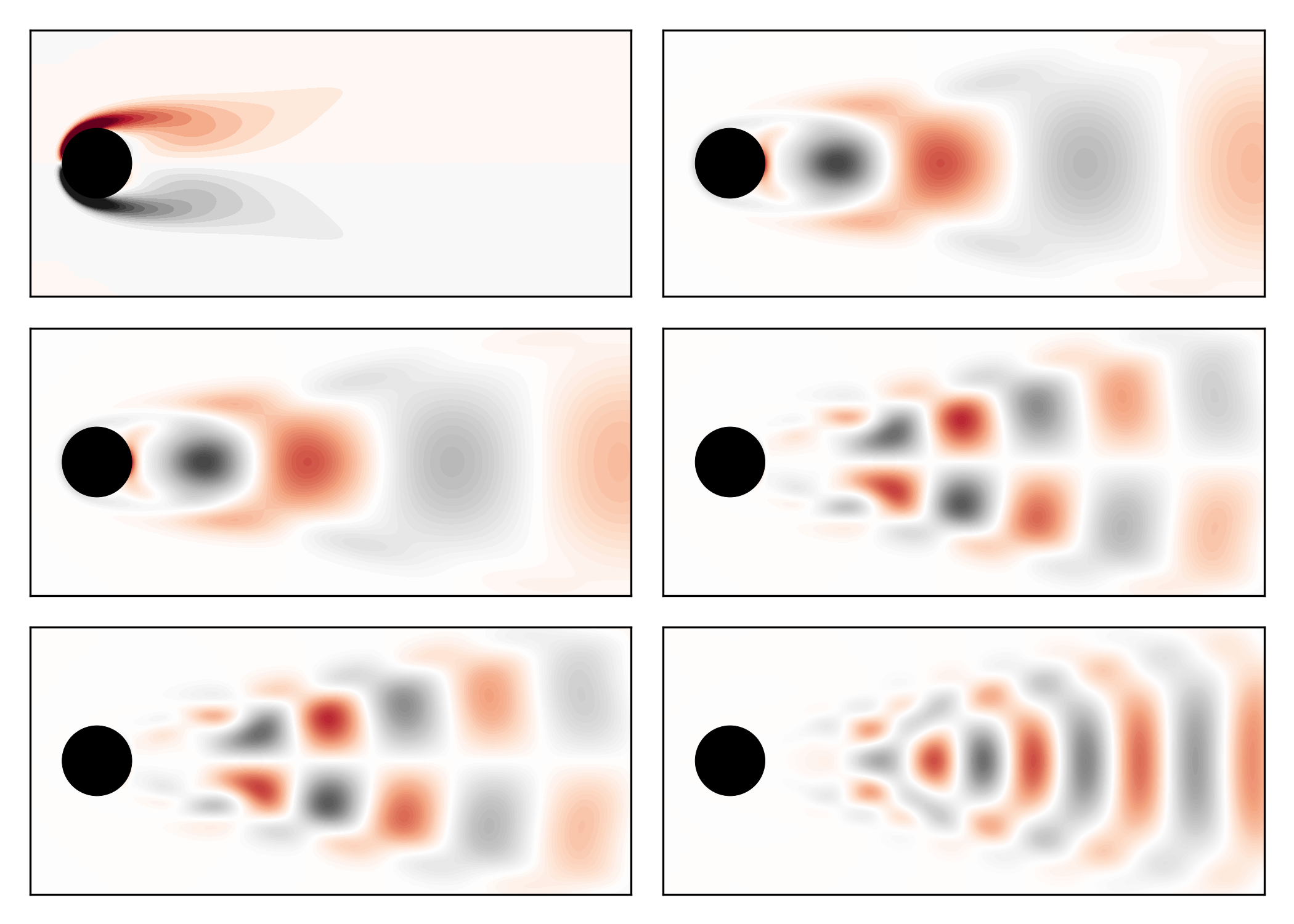}
		\caption{Deterministic DMD modes.}
	\end{subfigure}
	~
	\begin{subfigure}[t]{0.48\textwidth}
		\centering
		\DeclareGraphicsExtensions{.png}
		\includegraphics[width=1\textwidth]{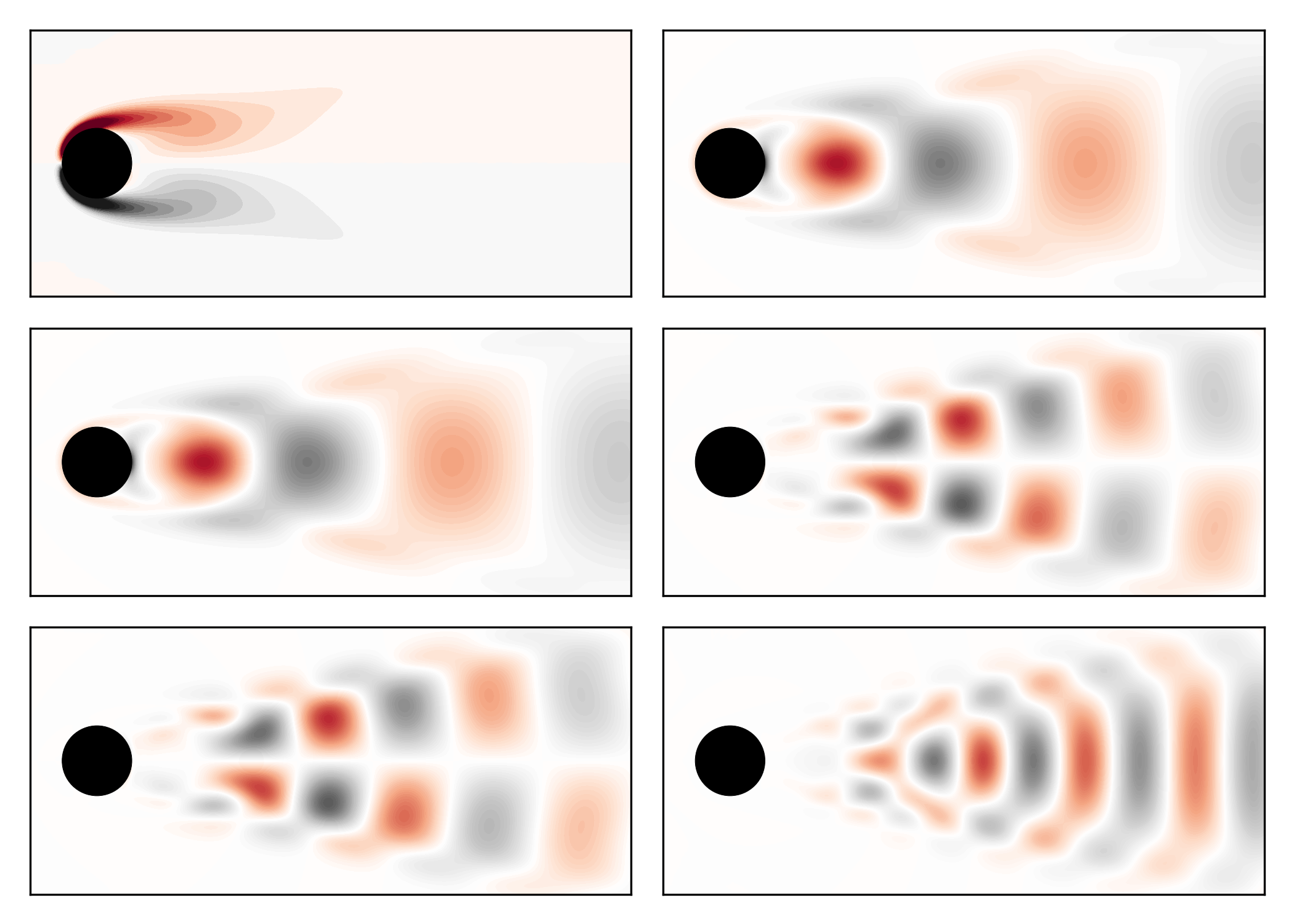}				
		\caption{Randomized DMD modes.}
	\end{subfigure}		
	
	\caption{Leading dynamic modes extracted from the fluid flow behind a cylinder.}
	\label{fig:flow_modes}
\end{figure}

We compute the low-rank DMD approximation using $k=15$ as the desired target rank.
Figure~\ref{fig:flow_eigA} shows the DMD eigenvalues. The proposed randomized DMD algorithm (with $p=10$ and $q=0$), and the compressed DMD algorithm (with $l=1000$) faithfully capture the eigenvalues.
Overall, the randomized algorithm leads to a $6$ fold speedup compared to the deterministic DMD algorithm.
Further, if the singular value spectrum is slowly decaying, the approximation accuracy can be improved by computing additional power iterations $q$.  
To further contextualize the results, Fig.~\ref{fig:flow_modes} shows the leading six DMD modes in absence of noise.
The randomized algorithm faithfully reveals the coherent structures, while requiring considerably fewer computational resources. 

%

Next, the analysis is repeated in presence of additive white noise with a signal-to-noise ratio (SNR) of $10$. Figure~\ref{fig:flow_eigB} shows distinct performance of the different algorithms. The deterministic algorithm performs most accurately, capturing the first eleven eigenvalues. The randomized DMD algorithm reliably captures the first nine eigenvalues, while the compressed algorithm only accurately captures seven of the eigenvalues.
This is, despite the fact that the compressed DMD algorithm uses a large sketched snapshot sequence of dimension $1000 \times 151$, i.e., $1,000$ randomly selected rows out of the $89,351$ rows of the flow data set are used.
For comparison the randomized DMD algorithm provides a more satisfactory approximation using $l=25$ and $2$ additional power iterations, i.e., the sketched snapshot sequence is only of dimension $25 \times 151$. 
The results show that the randomized DMD algorithm yields a more accurate approximation, while allowing for higher compression rates.
This is because the approximation quality of the randomized DMD algorithm depends on the intrinsic rank of the data and not on the ambient dimensions on the measurement space.

\begin{figure}[!b]
	\centering
	\begin{minipage}[b]{.48\textwidth}
		\begin{subfigure}[t]{0.99\textwidth}
			\centering
			\DeclareGraphicsExtensions{.png}
			\includegraphics[width=1\textwidth]{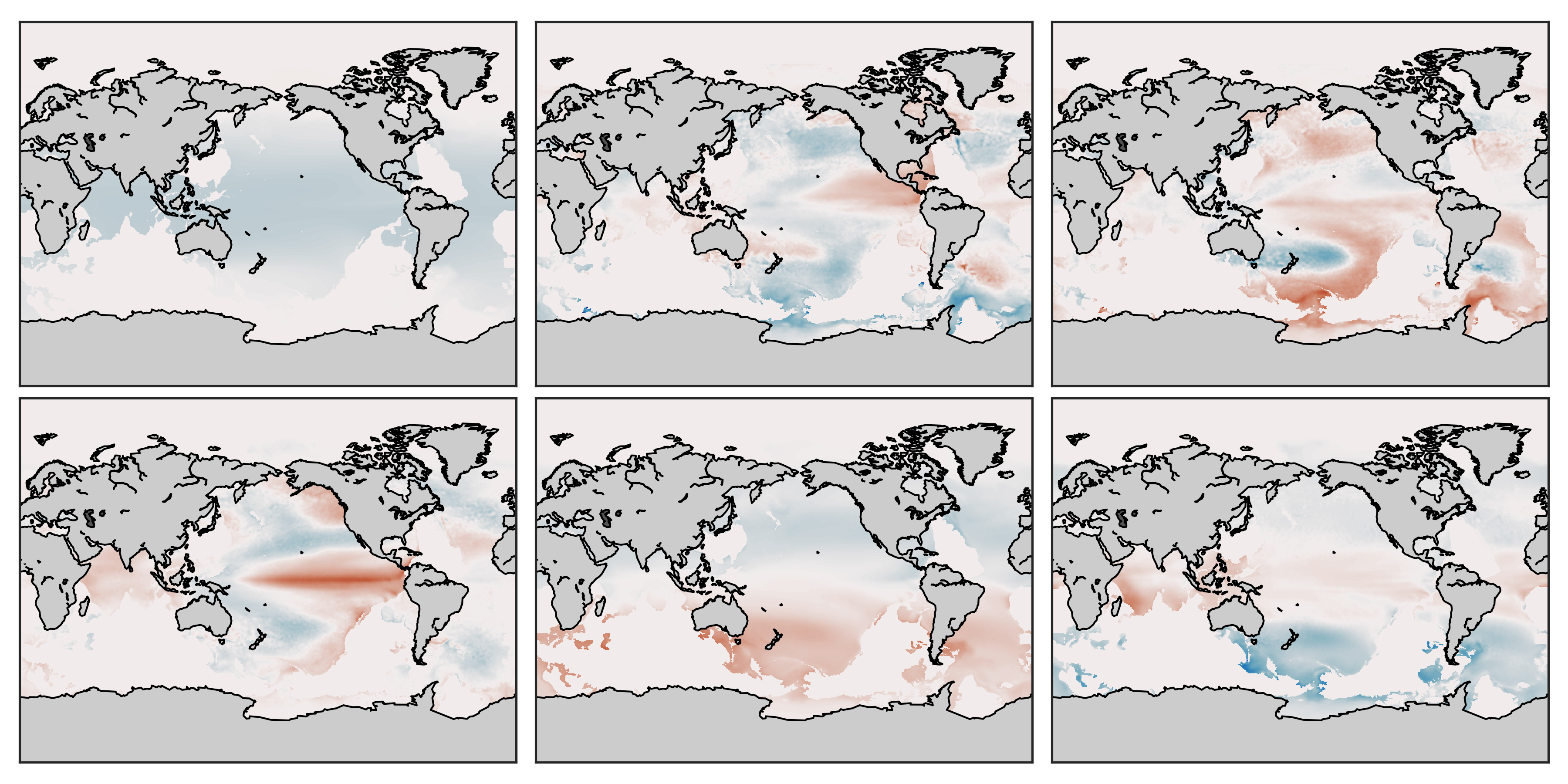}
			\caption{Deterministic DMD modes.}
			\label{fig:sst_modes}
			\includegraphics[width=1\textwidth]{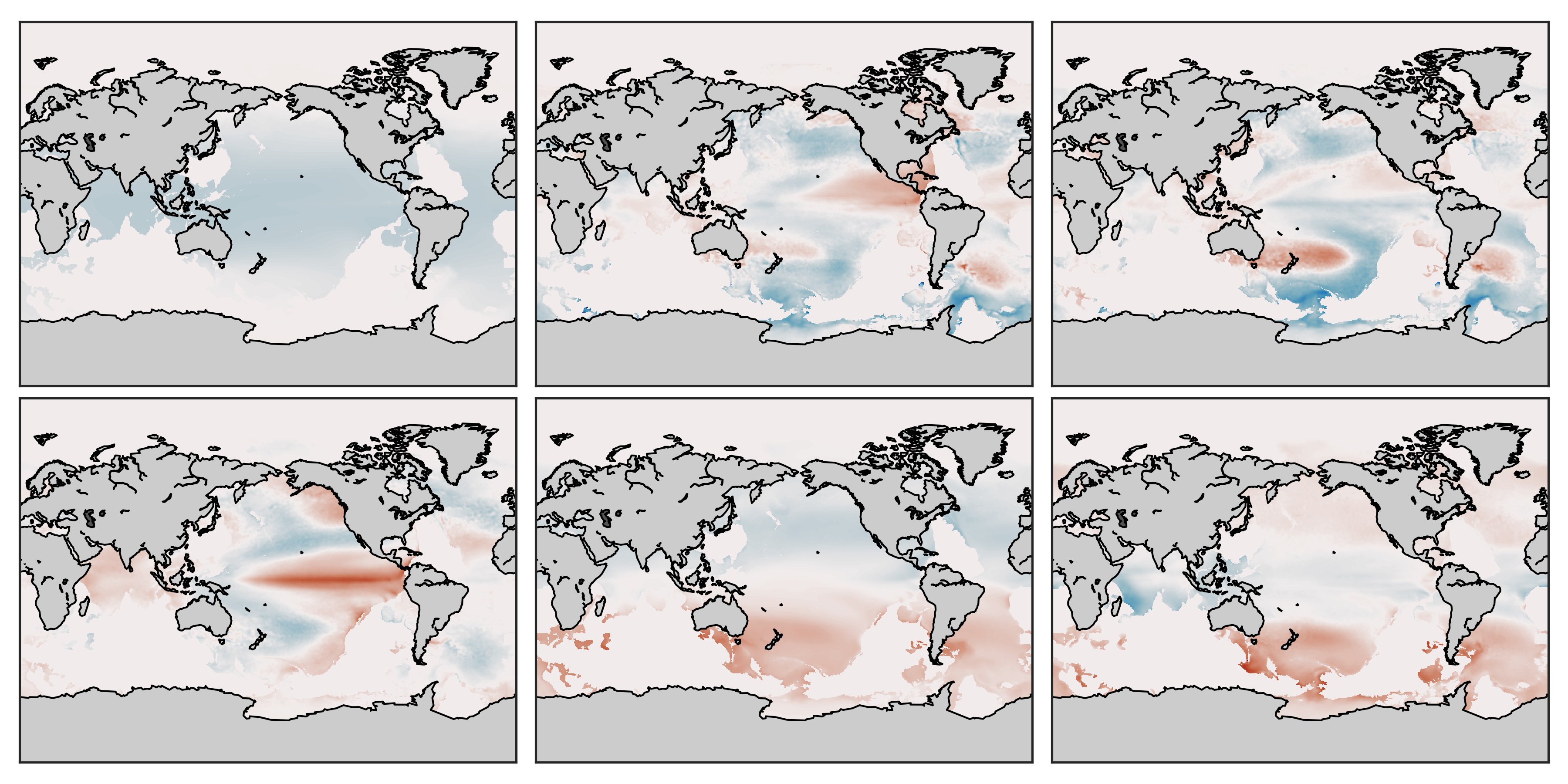}				
			\caption{Blocked Randomized DMD modes.}
			\label{fig:sst_modes_rand}
		\end{subfigure}
	\end{minipage}
	~
	\begin{minipage}[b]{.48\textwidth}
		\begin{subfigure}[b]{0.99\textwidth}
			\centering
			\DeclareGraphicsExtensions{.pdf}
			\includegraphics[width=1\textwidth]{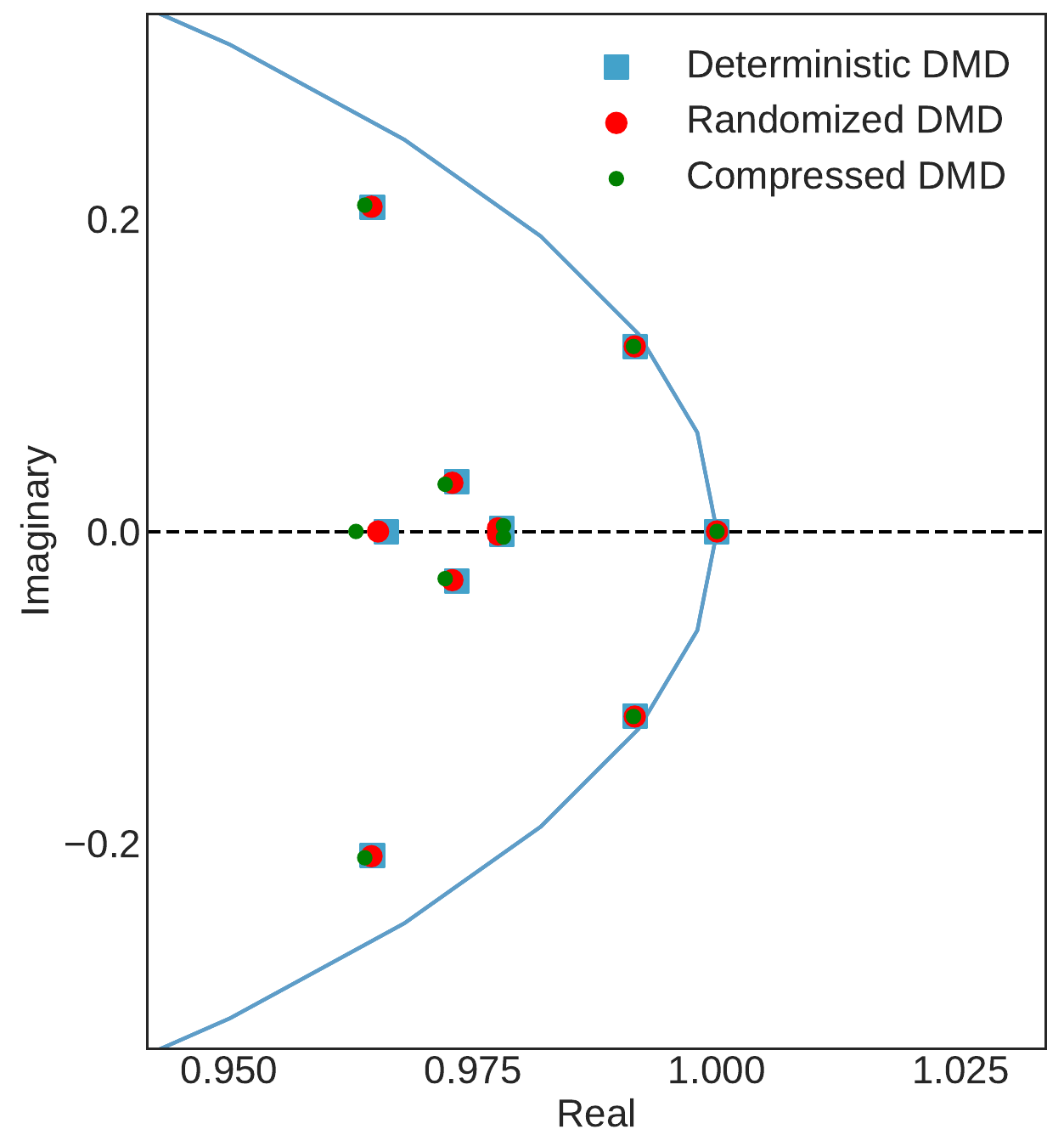}
			\caption{DMD eigenvalues.}
			\label{fig:fluid_eig}
		\end{subfigure}
		\vspace{-2.38in}
	\end{minipage}		
	\caption{Leading dynamic modes extracted from the aggregated high-resolution sea surface dataset. No distinct differences can be obtained. Further, the randomized DMD algorithm captures faithfully the eigenvalues. In contrast, the compressed DMD algorithm provides only a crude approximation for the eigenvalues. }
	\label{fig:sst_summary}
\end{figure}

\subsection{Sea Surface Data} 

We now compute the dynamic mode decomposition on the high-resolution sea surface temperature (SST) data. 
The SST data are widely studied in climate science for climate monitoring and prediction, providing an improved understanding of the interactions between the ocean and the atmosphere~\cite{reynolds2007daily,reynolds2002improved,smith2005global}.
Specifically, the daily SST measurements are constructed by combining infrared satellite data with observations provided by ships and buoys.  
In order to account and compensate for platform differences and sensor errors, a bias-adjusted methodology is used to combine the measurements from the different sources.
Finally, the spatially complete SST map is produced via interpolation. 
A comprehensive discussion of the data is provided in~\cite{reynolds2002improved}.

The data are provided by the National Oceanic and Atmospheric Administration (NOAA) via their web site at \url{https://www.esrl.noaa.gov/psd/}. Data are available for the years from $1981$ to $2018$ with a temporal resolution of 1 day and a spatial grid resolution of $0.25^\circ$.
In total, the data consist of $m=13,149$ temporal snapshots which measure the daily temperature at $1440\times 720 = 1,036,800$ spatial grid points. Since we omit data over land, the ambient dimension reduces to $n=691,150$ spatial measurements in our analysis. Concatenating the reduced data yield a $36$GB data matrix of dimension $\mathbf{X} \in \mathbb{R}^{691,150 \times 13,149}$, which is sufficiently large to test scaling.

\subsubsection{Aggregated (Weekly Mean) Data}

The full data set outstrips the available fast memory required to compute the deterministic DMD. Thus, we perform the analysis on an aggregated data set first in order to compare the randomized and deterministic DMD algorithms. Therefore, we compute the weekly mean temperature, which reduces the number of temporal snapshots to $m=1,878$. 
Figure~\ref{fig:fluid_eig} shows the corresponding eigenvalues. Unlike the eigenvalues obtained via the compressed DMD algorithm, the randomized DMD algorithm provides an accurate approximation for the dominant eigenvalues. 
Next, Fig.~\ref{fig:sst_modes} and~\ref{fig:sst_modes_rand} show the extracted DMD modes. Indeed, the randomized modes faithfully capture the dominant coherent structure in the data. 

%

\subsubsection{Full Data}

\begin{figure}[!b]
	\centering
	\DeclareGraphicsExtensions{.png}
	\includegraphics[width=.98\textwidth]{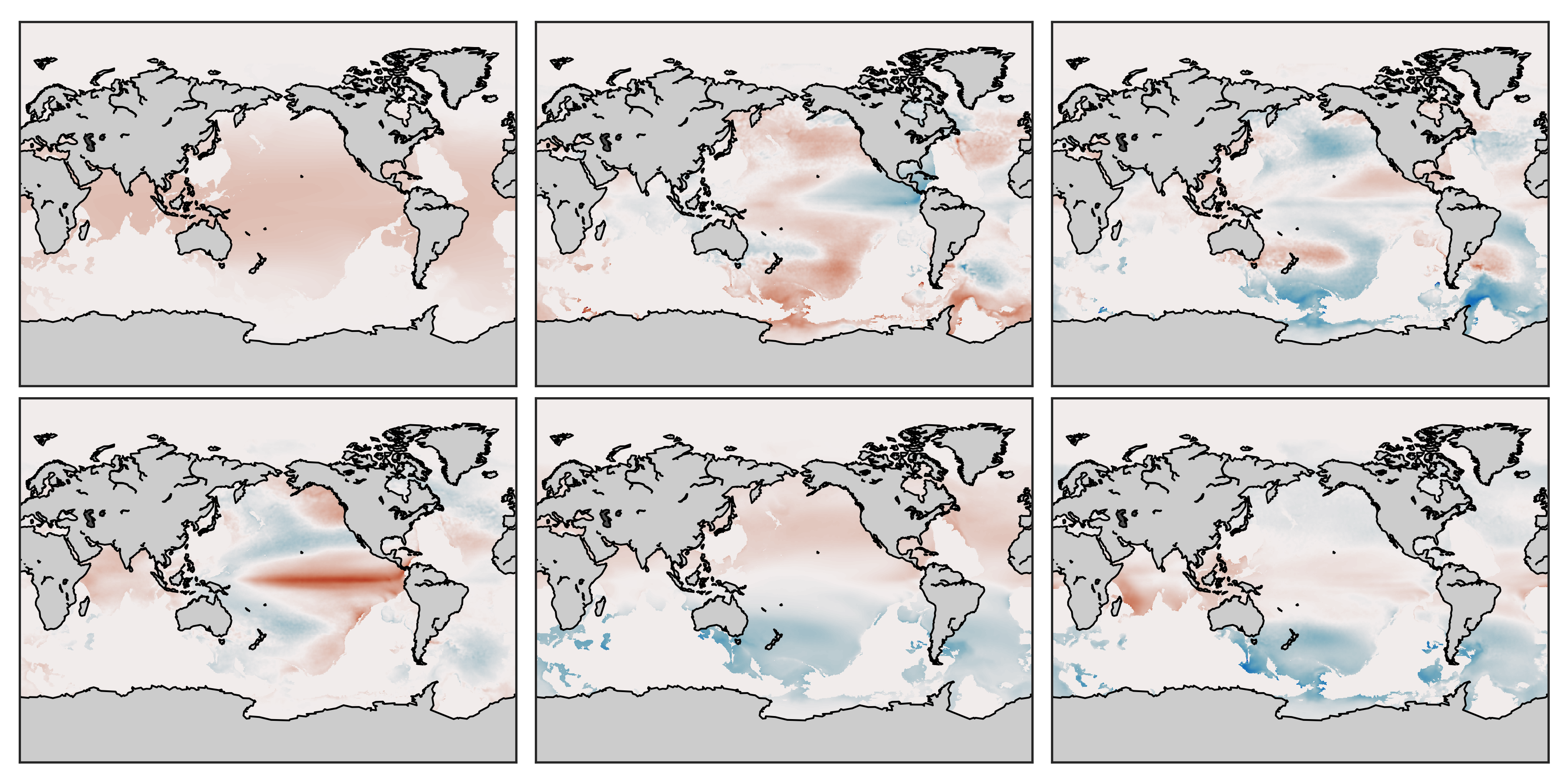}				
	\caption{Leading dynamic modes of the full high-resolution SST dataset.}
	\label{fig:sst_full_modes}
\end{figure}

The blocked randomized DMD algorithm allows us to extract the DMD modes from the high-dimensional data set. While the full data set does not fit into fast memory, it is only required that we can access some of the rows in a sequential manner. 
Computing the $k=15$ approximation using $b=4$ blocks takes about $150$ seconds.
The resulting modes are shown in Fig.~\ref{fig:sst_full_modes}. Here, the leading modes are similar to the modes extracted from the aggregated data set. However, subsequent modes may provide further insights which are not revealed by the aggregated data set.

\subsection{Computational Performance}

Next, we evaluate the computational performance of the randomized algorithms in terms of both time and accuracy. We measure the accuracy by computing the relative error of the randomized DMD compared to the approximation produced by the deterministic DMD algorithm
\begin{equation}
	\rho(\mathbf{\widehat{X}}_\mathrm{(DMD)}, \mathbf{\widehat{X}}_\mathrm{(rDMD)}) = \frac{\|\mathbf{\widehat{X}}_\mathrm{(DMD)} - \mathbf{\widehat{X}}_\mathrm{(rDMD)}  \|_F}{\|\mathbf{\widehat{X}}_\mathrm{(DMD)}\|_F},
\end{equation}
where $\mathbf{\widehat{X}}$ denotes the reconstructed snapshot matrix using either the deterministic or randomized DMD algorithm. Here, we use a (high-performance) partial SVD algorithm to compute the $k$ deterministic modes~\cite{lehoucq1997arpack}.
\begin{figure}[!b]
	\centering
	
	\begin{subfigure}[b]{0.99\textwidth}
		\centering
		\DeclareGraphicsExtensions{.pdf}
		\includegraphics[width=1\textwidth]{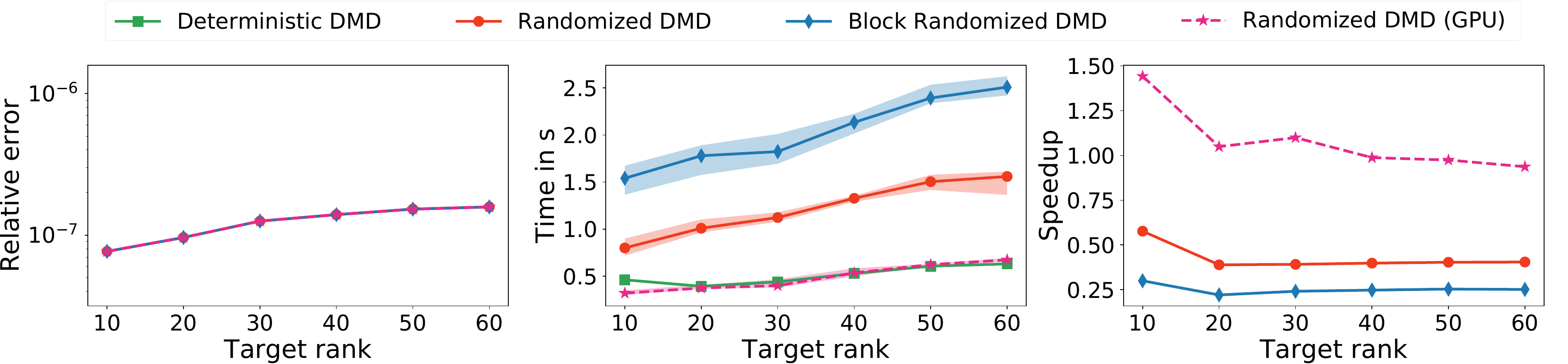}
		\caption{Flow past cylinder (matrix of dimension $89,351\times 151$).}
		\label{fig:flow_performance}
	\end{subfigure}	
	
	\begin{subfigure}[b]{0.99\textwidth}
		\centering
		\DeclareGraphicsExtensions{.pdf}
		\includegraphics[width=1\textwidth]{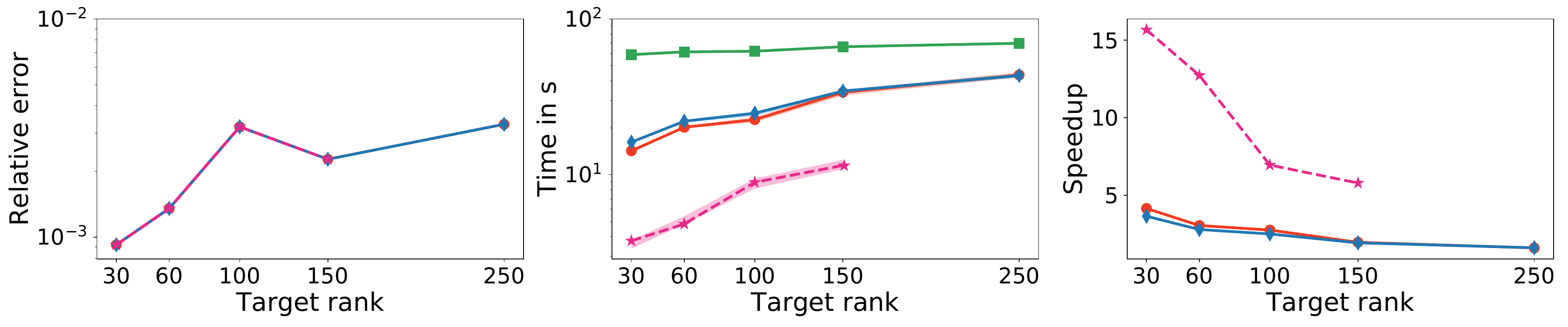}
		\caption{Aggregated sea surface temperature (matrix of dimension $691,150\times 1,878$).}
		\label{fig:sst_performance}
	\end{subfigure}		
	
	\begin{subfigure}[b]{0.99\textwidth}
		\centering
		\DeclareGraphicsExtensions{.pdf}
		\includegraphics[width=1\textwidth]{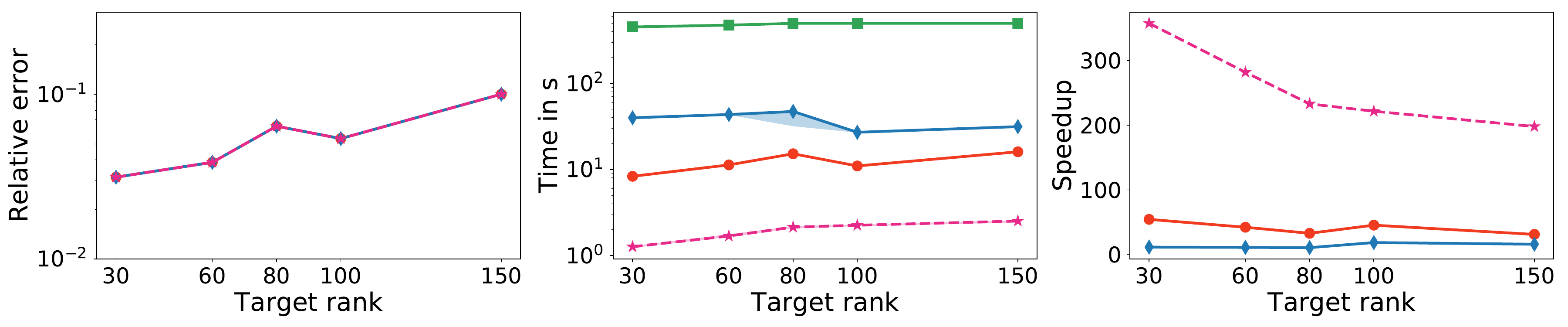}
		\caption{Turbulent flow (matrix of dimension $65,536\times 15,000$).}
		\label{fig:big_performance}
	\end{subfigure}	

	\begin{subfigure}[b]{0.99\textwidth}
	\centering
	\DeclareGraphicsExtensions{.pdf}
	\includegraphics[width=1\textwidth]{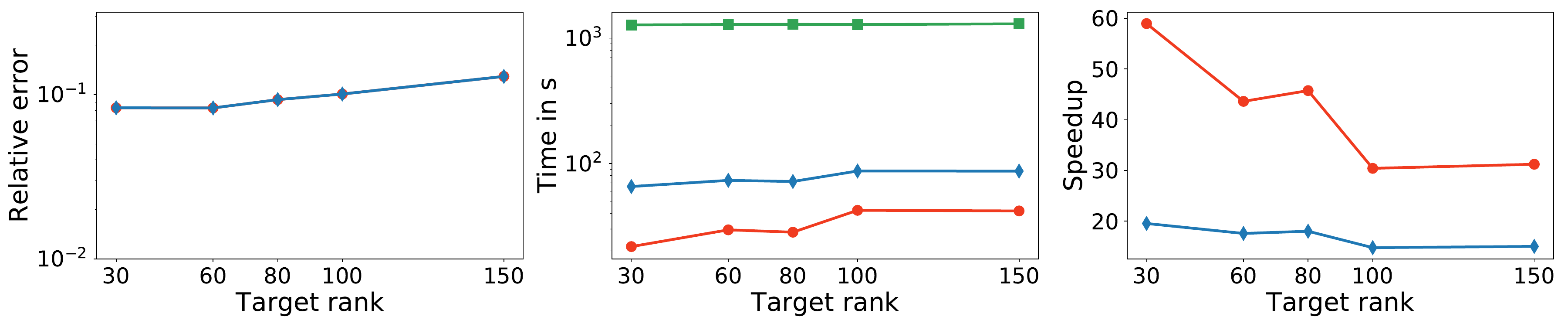}
	\caption{Turbulent flow (matrix of dimension $131,072\times 20,000$).}
	\label{fig:big_performance2}
	\end{subfigure}\vspace{-0.5cm}		
	
	\caption{Average runtimes and errors, over $20$ runs, for varying target ranks. Further, the gained speedup compared to the deterministic algorithm (baseline) is shown.}
	\label{fig:speed}
	
\end{figure}

\begin{figure}[!b]
	\centering
	\DeclareGraphicsExtensions{.pdf}
	\includegraphics[width=0.9\textwidth]{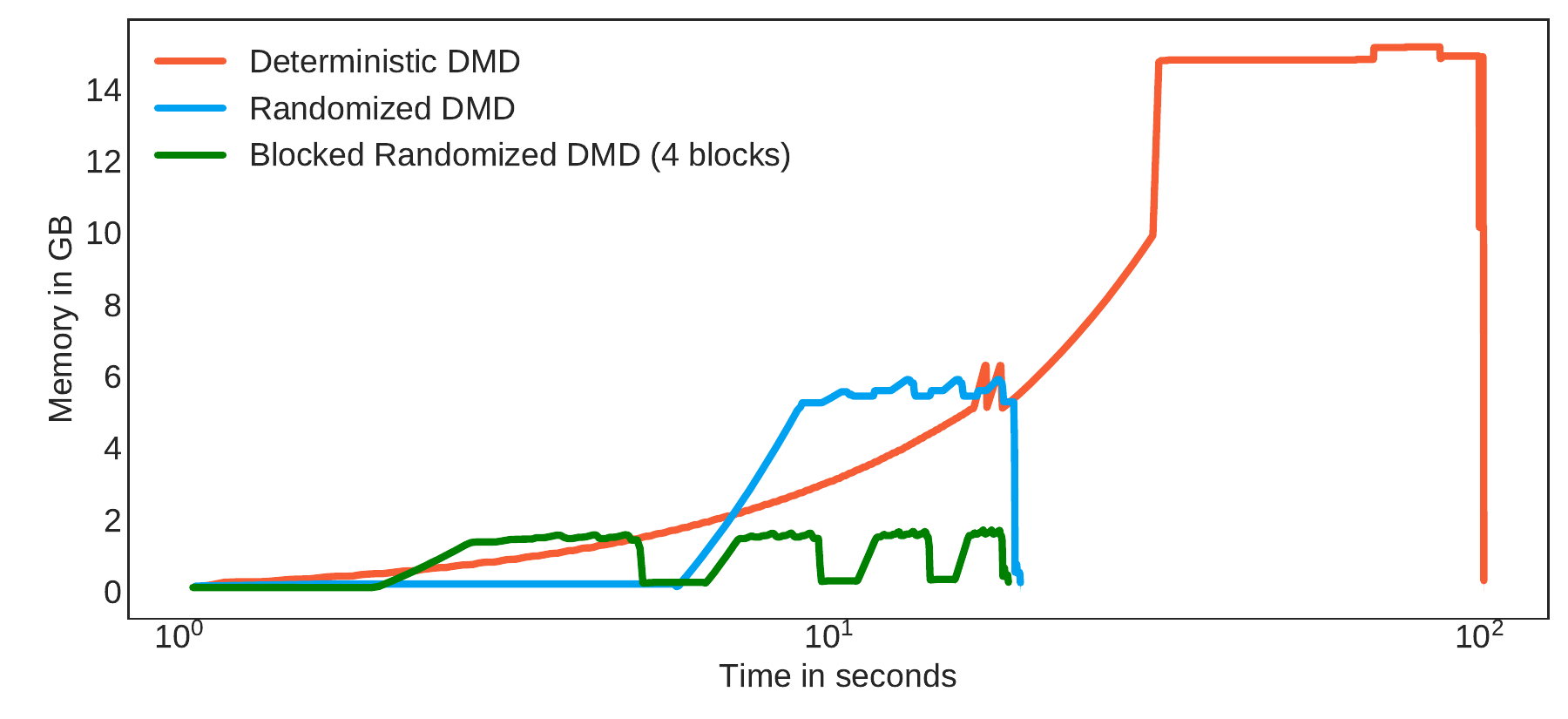}
	\caption{A profile of the memory usage vs runtime of the different DMD algorithms. }
	\label{fig:dmd_memory}
\end{figure}

Figure~\ref{fig:speed} summarizes the computational performance, for three different examples, and for varying target ranks.
First, Figure~\ref{fig:speed}~(a) shows the performance for the flow past a cylinder. The snapshot matrix is tall and skinny. In this setting there is no computational advantage of the randomized algorithm over the deterministic algorithm. This is because the data matrix is overall very small and the (deterministic) partial SVD is highly efficient for tall and skinny matrices. Importantly, we stress that the randomized DMD algorithm provides a very accurate approximation for the flow past cylinder. That is, because this fluid flow example has a fast decaying singular value spectrum.   

Second, Figure~\ref{fig:speed}~(b) shows the performance for the aggregated SST data. This snapshot matrix is also tall and skinny, yet its dimensions are substantially larger than those of the previous example. 
Here, we start to see the computational benefits of the randomized DMD algorithm. We gain about a speedup factor of 3-5 for computing the low-rank DMD, while maintaining a good accuracy. Note that this problem has a more slowly decaying singular value spectrum and thus the accuracy is poorer than in the previous examples. 
%
%
%
This example also demonstrates the performance boost by using a GPU accelerated implementation of the randomized DMD. Clearly, we can see a significant speedup by a factor of about 10-15 for computing the top 30 to 60 dynamic modes. Note that we run out of memory for target ranks larger than $k >150$.

Third, Figure~\ref{fig:speed}~(c)~and~(d) shows the performance for a turbulent flow data set. The fluid flow was obtained with the model implementation of~\cite{Chirila2018}, which is based on the algorithm presented by~\cite{Wang2013a}. This data set shows the advantages and disadvantages of the randomized algorithm. Indeed, we see some substantial gains in terms of the computational time, i.e., the GPU-accelerated algorithms achieve speedups of factors $>300$. However, the accuracy is less satisfactory due to the slowly decaying singular value spectrum of the data. This example illustrates the limits of the randomized algorithm, i.e., we need to assume that the data set features low-rank structure and has a reasonably fast decaying singular value spectrum. 

The reader might have noticed that the blocked randomized scheme (using 2 blocks) requires more computational time to compute the approximate DMD. This is because the data matrices considered in the above examples fit into the fast memory. Hence, there is little advantage in using the blocked scheme here.
However, the memory requirements are often more important than the absolute computational time. Figure~\ref{fig:dmd_memory} shows a profile of the memory usage vs runtime. The blocked randomized DMD algorithm requires only a fraction of the  memory, compared to the deterministic algorithm, to compute the approximate low-rank dynamic mode decomposition. This can be crucial when carrying out the computations on mobile platforms, on GPUs (as seen in Figure~\ref{fig:speed}~(c)), or when scaling the computations to very large applications.

\section[Conclusion]{Conclusion}\label{sec:conclusion}

Randomness as a computational strategy has recently been shown capable of efficiently solving many standard problems in linear algebra. The need for highly efficient algorithms becomes increasingly important in the area of `big data'. Here, we have proposed a novel randomized algorithm for computing the low-rank dynamic mode decomposition. Specifically, we have shown that DMD can be embedded in the probabilistic framework formulated by~\cite{halko2011rand}. This framework not only enables computations at scales far larger than what was previously possible, but it is also modular, flexible, and the error can be controlled. Hence, it can be also utilized as a framework to embed other innovations around the dynamic mode decomposition, for instance, see~\cite{Kutz2016book} for an overview. 

The numerical results show that randomized dynamic mode decomposition (rDMD) has computational advantages over previously suggested probabilistic algorithms for computing the dominant dynamic modes and eigenvalues. 
More importantly, we showed that the algorithm can be executed using a blocked scheme which is memory efficient. This aspect is crucial in order to efficiently deal with massive data which are too big to fit into fast memory. 
Thus, we believe that the randomized DMD framework will provide a powerful and scalable architecture for extracting dominant spatiotemporal coherent structures and dynamics from increasingly large-scale data, for example from epidemiology, neuroscience, and fluid mechanics.


\section*{Acknowledgments}
We would like to express our gratitude to the two anonymous reviewers. Their helpful feedback allowed us to greatly improve the manuscript. Further, we thank Diana A. Bistrian for many useful comments.
NBE would like to acknowledge the generous funding support from the Defense Advanced Research Projects Agency (DARPA) and the Air Force Research Laboratory (FA8750-17-2-0122) as well as Amazon Web Services for supporting the project with EC2 credits.
LM acknowledges the support of the French Agence Nationale pour la Recherche (ANR) and Direction G\'en\'erale de l'Armement (DGA) via the \textit{FlowCon} project (ANR-17-ASTR-0022).
SLB and JNK would like to acknowledge generous funding support from the Army Research Office (ARO W911NF-17-1-0118), Air Force Office of Scientific Research (AFOSR FA9550-16-1-0650), and DARPA (PA-18-01-FP-125). 

\bibliographystyle{siamplain}
\bibliography{references}   

\end{document}

%% file: ex_shared.tex
\usepackage{lipsum}
\usepackage{amsfonts}
\usepackage{graphicx}
\usepackage{epstopdf}
\usepackage{algorithmic}

\usepackage{url}            
\usepackage{amsfonts}       
\usepackage{microtype}      

\usepackage{amsmath}
\usepackage{amssymb}
\usepackage{graphicx}
\usepackage{epstopdf}
\usepackage{subcaption}
\usepackage{multirow}
\usepackage{float}
\usepackage{setspace}
\usepackage{overpic}
\usepackage{rotating}

\usepackage{tikz}
\usetikzlibrary{matrix,positioning,shapes,shadows,arrows,calc,decorations.pathreplacing}

\tikzstyle{block} = [draw,line width=1.2pt, fill=white!30, rectangle, 
minimum height=3em, minimum width=3em, rounded corners]         
\tikzstyle{blocky} = [draw, line width=1.2pt,fill=yellow!40, rectangle, 
minimum height=2.5em, minimum width=2.5em, rounded corners]    
\tikzstyle{blockp} = [draw,line width=1.2pt, fill=purple!30, rectangle, 
minimum height=3em, minimum width=3em, rounded corners]       
\tikzstyle{blockb} = [draw,line width=1.2pt, fill=blue!20, rectangle, 
minimum height=2.5em, minimum width=2.5em, rounded corners]     
\tikzstyle{blockr} = [draw, line width=1.2pt,fill=red!30, rectangle, 
minimum height=2.5em, minimum width=2.5em, rounded corners]                 
\tikzstyle{blockgr} = [draw, line width=1.2pt,fill=black!10, rectangle, 
minimum height=2.5em, minimum width=2.5em, rounded corners]                  
\tikzstyle{line}=[-, line width=1.2pt]

\tikzstyle{blocky4} = [draw, line width=1.2pt,fill=white!10, rectangle, 
text centered,minimum height=5.0em, text width=6.5em, rounded corners]        

\tikzstyle{blocky5} = [draw, line width=0pt, fill=white!0, opacity=0.0, rectangle, 
text centered,minimum height=5.0em, text width=6.5em, rounded corners]   

\tikzstyle{blockyY} = [draw, line width=1.2pt,fill=red!10, rectangle, 
text centered,minimum height=2.5em, text width=6.5em, rounded corners]    

\tikzstyle{blockyZ} = [draw, line width=1.2pt,fill=red!10, rectangle, 
text centered,minimum height=2.5em, text width=4.0em, rounded corners] 

\tikzstyle{blockyC} = [draw, line width=1.2pt,fill=red!10, rectangle, 
text centered,minimum height=2.5em, text width=4.0em, rounded corners] 

\DeclareMathOperator*{\argmin}	{{\mathrm{ arg\,min}}}

\usepackage{booktabs}

\newcommand\Rey{\mbox{\textit{Re}}}  

\newcommand{\bx}{\mathbf{x}}

\newcommand{\bA}{\mathbf{A}}

\newcommand{\bP}{\mathbf{P}}
\newcommand{\bT}{\mathbf{T}}
\newcommand{\bU}{\mathbf{U}}
\newcommand{\bV}{\mathbf{V}}
\newcommand{\bW}{\mathbf{W}}
\newcommand{\bX}{\mathbf{X}}
\newcommand{\bY}{\mathbf{Y}}

\newcommand{\bQ}{\mathbf{Q}}
\newcommand{\bB}{\mathbf{B}}
\newcommand{\bI}{\mathbf{I}}
\newcommand{\bb}{\mathbf{b}}

\newcommand{\bSig}{\mathbf{\Sigma}}
\newcommand{\bK}{\mathbf{K}}

\newcommand{\bSigma}{\boldsymbol{\Sigma}}

\DeclareMathOperator{\E}{\mathbb{E}}

\ifpdf
  \DeclareGraphicsExtensions{.eps,.pdf,.png,.jpg}
\else
  \DeclareGraphicsExtensions{.eps}
\fi

\usepackage{enumitem}
\setlist[enumerate]{leftmargin=.5in}
\setlist[itemize]{leftmargin=.5in}

\newsiamremark{remark}{Remark}
\newsiamremark{hypothesis}{Hypothesis}
\crefname{hypothesis}{Hypothesis}{Hypotheses}
\newsiamthm{claim}{Claim}

\headers{Randomized Dynamic Mode Decomposition}{N. Benjamin Erichson, Lionel Mathelin, J. Nathan Kutz, Steven L. Brunton}

\title{Randomized Dynamic Mode Decomposition}

\author{N. Benjamin Erichson\footnotemark[2] \footnotemark[5]
\and Lionel Mathelin\footnotemark[3]\ \footnotemark[2]
\and J. Nathan Kutz\footnotemark[2]
\and Steven L. Brunton\footnotemark[4]}

\footnotetext[2]{Department of Applied Mathematics, University of Washington, Seattle, WA ({kutz@uw.edu})}
\footnotetext[5]{ICSI and Department of Statistics, University of California, Berkeley ({erichson@berkeley.edu})}
\footnotetext[3]{LIMSI--CNRS, Campus Universitaire d'Orsay, 91405 Orsay cedex, France ({mathelin@limsi.fr})}
\footnotetext[4]{Department of Mechanical Engineering, University of Washington, Seattle, WA ({sbrunton@uw.edu})}

\usepackage{amsopn}
\DeclareMathOperator{\diag}{diag}
